\documentclass[a4paper]{article}%
\usepackage{amsmath}
\usepackage{graphicx}
\usepackage{amsfonts}
\usepackage{amssymb}%
\begin{document}

\title{Schubert calculus in Lie groups}
\author{Haibao Duan}
\date{}
\maketitle

\begin{abstract}
Let $G$ be a Lie group with a maximal torus $T$. Combining Schubert calculus
in the flag manifold $G/T$ with the Serre spectral sequence of the fibration
$G\rightarrow G/T$, we construct the integral cohomology ring $H^{\ast}(G)$
uniformly for any compact and simply-connected Lie group $G$.

\begin{description}
\item[2000 Mathematical Subject Classification: ] 14M15; 55T10; 57T10

\item[Key words and phrases:] Lie groups; Cohomology, Schubert calculus, Serre
spectral sequence, Koszul complex

\end{description}
\end{abstract}

\section{Introduction}

The problem of computing the cohomology of Lie groups was raised by \'{E}.
Cartan in 1929 \cite{C,K}. It is a focus of algebraic topology owing to the
fundamental role of Lie groups playing in geometry and topology \cite{D,MT,Sa}%
. However, despite contributions by many mathematicians in about one century,
notably \cite{B2,B3,BC,Br,G,H,K,M,M2,Po}, the integral cohomology of Lie
groups remains incomplete. On the other hand, Schubert calculus began with the
enumerative geometry of the 19th century. Making this calculus rigorous was
stated by Hilbert as his 15th problem. In the course of securing the
foundation of algebraic geometry, Van der Waerden \cite{Wa} and A.Weil
\cite[p.331]{W} related the problem to the determination of the intersection
theory of flag manifolds, see \cite{DZ3,DZ4} for surveys about the history. In
this paper we bring a connection between these two topics both with
distinguished backgrounds, and illustrate how Schubert calculus in flag
manifolds can be extended to obtain an explicit and unified construction of
the integral cohomology rings of all compact and simply-connected Lie groups,
e.g. Theorem B in Section \S 1.3.

It is well known that any compact and simply-connected Lie group $G$ is
isomorphic to a product $G_{1}\times\cdots\times G_{k}$, where each factor
$G_{i}$ is a simply-connected and simple Lie group. It is also known that all
the simply-connected and simple Lie groups fall into three infinite families
of the classical Lie groups $SU(n)$, $Spin(n),Sp(n)$, as well as the five
exceptional ones $G_{2},F_{4},E_{6},E_{7},E_{8}$. Therefore, we shall assume
in this paper that $G$ is one of these simple Lie group. The cohomologies are
over the ring $\mathbb{Z}$ of integers, unless otherwise stated.

\subsection{Preliminaries in Schubert calculus}

Fixing a maximal torus $T$ on $G$ the inclusion $T\subset G$ induces the fiber sequence

\begin{enumerate}
\item[(1.1)] $T\rightarrow G\overset{\pi}{\rightarrow}G/T\overset
{c}{\rightarrow}B_{T}$,
\end{enumerate}

\noindent where $B_{T}$ is the classifying space of the group $T$, and the
quotient space $G/T$ is a complex projective manifold, called \textsl{the
complete flag manifold} of $G$ \cite{DZ2,DZ3}. Suppose that the rank of $G$ is
$\dim T=n$, and $\{\omega_{1},\cdots,\omega_{n}\}\subset H^{2}(B_{T})$ is a
set of fundamental dominant weights of the group $G$ \cite{BH}. Then the
cohomology $H^{\ast}(B_{T})$ is isomorphic to the polynomial ring
$\mathbb{Z}[\omega_{1},\cdots,\omega_{n}]$, while the induced map of $c$ is
the \textsl{Borel characteristic map} of $G$ \cite{De1,De2}

\begin{enumerate}
\item[(1.2)] $c^{\ast}:\mathbb{Z}[\omega_{1},\cdots,\omega_{n}]\rightarrow
H^{\ast}(G/T)$.
\end{enumerate}

\noindent Furthermore, since the ring $H^{\ast}(G/T)$ is torsion free
\cite{BS}, the second page of the Serre spectral sequence $\{E_{r}^{\ast,\ast
}(G),d_{r}\}$ of $\pi$ has the presentation

\begin{enumerate}
\item[(1.3)] $E_{2}^{\ast,\ast}(G)=H^{\ast}(G/T)\otimes H^{\ast}(T)$,
\end{enumerate}

\noindent on which the $d_{2}$-action has been determined in \cite[Corollary
3.1]{D1}.

\bigskip

\noindent\textbf{Lemma 1.1. }\textsl{There exists a basis }$\left\{
\lambda_{1},\cdots,\lambda_{n}\right\}  $\textsl{\ of }$H^{1}(T)$%
\textsl{\ such that }

\textsl{i)} $H^{\ast}(T)=$\textsl{\ }$\Lambda(\lambda_{1},\cdots,\lambda_{n}%
)$\textsl{;}

\textsl{ii) }$d_{2}(x\otimes1)=0$\textsl{, }$d_{2}(x\otimes\lambda
_{i})=c^{\ast}(\omega_{i})\cdot x\otimes1$\textsl{, }$1\leq i\leq n$\textsl{,}

\noindent\textsl{where }$x\in H^{\ast}(G/T)$\textsl{,} \textsl{and }%
$\Lambda(\lambda_{1},\cdots,\lambda_{n})$\textsl{ is the exterior algebra
generated by }$\lambda_{1},\cdots,\lambda_{n}.$\hfill$\square$\noindent

\bigskip

Our aim is to construct the cohomology $H^{\ast}(G)$ in term of $E_{2}%
^{\ast,\ast}(G)$. To this end a concise characterization of the ring $H^{\ast
}(G/T)$ is requested. Early in the 1950, Borel and Chevalley computed the
cohomology $H^{\ast}(G/T;\mathbb{R})$ with real coefficients. Precisely,
assume that the Weyl group of $G$ is $W\subset Aut(H^{2}(B_{T}))$, and let
$\{q_{1},\cdots,q_{n}\}\in\mathbb{R}[\omega_{1},\cdots,\omega_{n}]^{W}$ be a
set of basic homogeneous $W$-invariants of $G$ \cite{K}, to be arranged in the
order $l_{1}\leq\cdots\leq l_{n}$, where $l_{i}$ is the degree of $q_{i}$ as a
polynomial in $\omega_{1},\cdots,\omega_{n}$. Borel and Chevalley \cite{B1,C}
showed that

\bigskip

\noindent\textbf{Lemma 1.2.} \textsl{The map }$c^{\ast}$ \textsl{induces an
isomorphism of algebras}

\begin{enumerate}
\item[(1.4)] $H^{\ast}(G/T;\mathbb{R})=\mathbb{R}[\omega_{1},\cdots,\omega
_{n}]/\left\langle q_{1},\cdots,q_{n}\right\rangle $\textsl{,}
\end{enumerate}

\noindent\textsl{where} \textsl{the} \textsl{set }$q(G)=\{l_{1},\cdots
,l_{n}\}$ \textsl{is an invariant of }$G$\textsl{, which is listed below:}

\begin{center}
{\normalsize Table 1. The degree set }$q(G)$ {\normalsize of the simple Lie
groups }${\normalsize G}$
\end{center}

\begin{quote}%
\begin{tabular}
[c]{l||l}\hline
$G$ & $q(G)=\{l_{1},\cdots,l_{n}\}$\\\hline\hline
$SU(n)$ & $\{2,3,\cdots,n\}$\\\hline
$Sp(n),Spin(2n+1)$ & $\{2,4,\cdots,2n\}$\\\hline
$Spin(2n)$ & $\{2,4,\cdots,2n-2,n\}$\\\hline
$G_{2}$ & $\{2,6\}$\\\hline
$F_{4}$ & $\{2,6,8,12\}$\\\hline
$E_{6}$ & $\{2,5,6,8,9,12\}$\\\hline
$E_{7}$ & $\{2,6,8,10,12,14,18\}$\\\hline
$E_{8}$ & $\{2,8,12,14,18,20,24,30\}$\\\hline
\end{tabular}
$.$\hfill$\square$\noindent
\end{quote}

Turning to the integral cohomology $H^{\ast}(G/T)$ we resorts to the classical
Schubert calculus. According to Chevalley \cite{Ch2}, or
Bernstein-Gel'fand-Gel'fand \cite{BGG}, the flag manifold $G/T$ admits a
decomposition into the Schubert cells $S_{w}$, parameterized by the elements
of the Weyl group $W$,

\begin{enumerate}
\item[(1.5)] $G/T=\cup_{w\in W}S_{w},\quad\dim S_{w}=2\iota(w)$,
\end{enumerate}

\noindent where $\iota:W\rightarrow\mathbb{Z}$ is the length function on $W$
\cite{BGG}. Since only even dimensional cells are involved, the set
$\{[S_{w}],w\in W\}$ of fundamental classes forms a basis of the homology
$H_{\ast}(G/T)$. The co-cycle class $s_{w}\in H^{\ast}(G/T)$ Kronecker dual to
the basis element $[S_{w}]$ (i.e. $\left\langle s_{w},[S_{u}]\right\rangle
=\delta_{w,u},~w,u\in W$) is called the \textsl{Schubert class} associated to
$w\in W$. Formula (1.5) implies the following result, well-known as "the basis
theorem of Schubert calculus" \cite{BGG,Ch2}.

\bigskip

\noindent\textbf{Lemma 1.3.} \textsl{The cohomology }$H^{\ast}(G/T)$
\textsl{is torsion free, concentrated in even degrees, and has} \textsl{a
basis\ consisting of the Schubert classes }$\{s_{w},w\in W\}$\textsl{.}

\textsl{Moreover, the Schubert basis of} $H^{2}(G/T)$ \textsl{is} $\{c^{\ast
}(\omega_{1}),\cdots,c^{\ast}(\omega_{n})\}$\textsl{\ (\cite{D3}).}%
\hfill$\square$\noindent

\bigskip

By the second assertion of Lemma 1.3 we may reserve the notion $\omega_{i}$
for $c^{\ast}(\omega_{i})\in H^{2}(G/T)$, and write $\left\langle \omega
_{1},\cdots,\omega_{n}\right\rangle _{c}$ to denote the ideal of $H^{\ast
}(G/T)$ generated by the $\omega_{i}$'s. Carrying on the results of Lemma 1.3
we show that

\bigskip

\noindent\textbf{Lemma 1.4.} \textsl{For each Lie group} $G$ \textsl{there
exists a minimal set of Schubert classes }$y_{t_{1}},\cdots,y_{t_{k}}$
\textsl{on} $G/T$\textsl{ with} $2<\deg y_{t_{1}}\leq\cdots\leq\deg y_{t_{k}}%
$\textsl{,} \textsl{which, together with} \textsl{the weights} $\omega
_{1},\cdots,\omega_{n}$\textsl{, generates} \textsl{the ring} $H^{\ast}%
(G/T)$\textsl{\ multiplicatively.}

\textsl{In addition, for each }$y_{t_{i}}$ \textsl{there exists a pair of}
\textsl{integers }$p_{i}$\textsl{,} $r_{i}>1$\textsl{, such that the following
relations hold in }$H^{\ast}(G/T)$

\begin{enumerate}
\item[(1.6)] $p_{i}y_{t_{i}}$\textsl{,} $y_{t_{i}}^{r_{i}}\in\left\langle
\omega_{1},\cdots,\omega_{n}\right\rangle _{c}$\textsl{.}
\end{enumerate}

\noindent\textbf{Proof. }Let $H^{+}(G/T)\subset H^{\ast}(G/T)$ be the subring
consisting of the homogeneous elements in positive degrees, and let
$D(H^{\ast}(G/T)):=H^{+}(G/T)\cdot H^{+}(G/T)$ be the ideal of the
decomposable elements of $H^{+}(G/T)$. Since $H^{\ast}(G/T)$ has a basis
consisting of Schubert classes, there exist Schubert classes $s_{1}%
,\cdots,s_{m}$, ordered by $\deg s_{1}\leq\cdots\leq$ $\deg s_{m}$, that
correspond to a basis of the quotient group $H^{+}(G/T)/D(H^{\ast}(G/T))$. It
follows that $\left\{  s_{1},\cdots,s_{m}\right\}  $ is a minimal set of
generators of the ring $H^{\ast}(G/T)$, and that the number $m$ is an
invariant of $G$.

By Lemma 1.3 the elements of $D(H^{\ast}(G/T))$ have degrees $\geq4$. In particular

\begin{quote}
$H^{2}(G/T)\subset H^{+}(G/T)/D(H^{\ast}(G/T))$.
\end{quote}

\noindent This implies by the second assertion of Lemma 1.3 that $m\geq n$,
and that we can take $s_{i}=\omega_{i}$, $1\leq i\leq n$. The proof of the
first assertion is completed by taking $y_{t_{i}}:=s_{n+i}$, where $1\leq
i\leq k:=m-n$.

By Lemma 1.2 the map $c^{\ast}$ surjects over $\mathbb{R}$. In particular, for
any $w\in W$ there exist positive integers $p_{w}$ and $r_{w}$ such that

\begin{quote}
$p_{w}s_{w}$, $s_{w}^{r_{w}}\in\left\langle \omega_{1},\cdots,\omega
_{n}\right\rangle _{c}$($\subseteq D(H^{\ast}(G/T))$).
\end{quote}

\noindent It implies, in particular, that if either $p_{w}=1$ or $r_{w}=1$, then

\begin{quote}
$s_{w}\equiv0\operatorname{mod}D(H^{\ast}(G/T))$.
\end{quote}

\noindent We obtain in (1.6) that $p_{i}$\textsl{,} $r_{i}>1$ from $y_{t_{i}%
}\neq0\operatorname{mod}D(H^{\ast}(G/T))$.\hfill$\square$\noindent

\bigskip

\noindent\textbf{Definition 1.5.} A set of Schubert classes $\left\{
y_{t_{1}},\cdots,y_{t_{k}}\right\}  $ satisfying Lemma 1.4 will be called
\textsl{a set of special Schubert classes} on $G/T$ (e.g. \cite{DZ2})

For each special Schubert class $y_{t_{i}}$ denote by $tor(y_{t_{i}})$ and
$cl(y_{t_{i}})$ the least multiple $p_{i}$ and the least power $r_{i}$
satisfying the relations (1.6), to be called \textsl{the torsion index} and
\textsl{the cup length} of $y_{t_{i}}$, respectively.\hfill$\square$\noindent

\bigskip

\noindent\textbf{Example 1.6. }By the algorithm indicated by the proof of
Lemma 1.4, and in the context of Schubert calculus \cite{D0}, we have selected
in \cite{DZ2} for every Lie group $G$ a set of special Schubert classes
$\left\{  y_{t_{1}},\cdots,y_{t_{k}}\right\}  $ on $G/T$, whose degree set
$\left\{  \deg y_{t_{1}},\cdots,\deg y_{t_{k}}\right\}  $, as well as the
indices $tor(y_{t_{i}})$ and lengths $cl(y_{t_{i}})$, are recorded in Table 2
below. We shall see in Theorem 4.1 that, these numbers are independent of the
choice of a set of special Schubert classes on $G/T$, hence are in fact
invariants of the group $G$, in addition to the degree set $q(G)$.

\begin{center}
{\normalsize Table 2. The degrees, torsion indices and cup lengths of the
special Schubert classes }$\left\{  y_{t_{1}},\cdots,y_{t_{k}}\right\}  $
{\normalsize on }${\normalsize G/T}$

{\footnotesize
\begin{tabular}
[c]{l||l|l|l|l}\hline
$G$ & $SU(n)$ & $Sp(n)$ & $Spin(2n)$ & $Spin(2n+1)$\\\hline\hline
$k$ & $0$ & $0$ & $[\frac{n-2}{2}]$ & $[\frac{n-1}{2}]$\\\hline
$\left\{  \deg y_{t_{1}},\cdots,\deg y_{t_{k}}\right\}  $ &  &  &
$\{4i+2,1\leq i\leq\lbrack\frac{n-2}{2}]\}$ & $\{4i+2,1\leq i\leq\lbrack
\frac{n-1}{2}]\}$\\\hline
$\{tor(y_{t_{i}})\}$ &  &  & $\{2,\cdots,2\}$ & $\{2,\cdots,2\}$\\\hline
$\{cl(y_{t_{i}})\}$ &  &  & $\{2^{[\ln\frac{n-1}{2i+1}]+1},1\leq i\leq
\lbrack\frac{n-2}{2}]\}$ & $\{2^{[\ln\frac{n}{2i+1}]+1},1\leq i\leq
\lbrack\frac{n-1}{2}]\}$\\\hline
\end{tabular}
}

\bigskip

{\footnotesize
\begin{tabular}
[c]{l}\hline%
\begin{tabular}
[c]{l||l|l|l|l|l}%
$G$ & $G_{2}$ & $F_{4}$ & $E_{6}$ & $E_{7}$ & $E_{8}$\\\hline\hline
$k$ & $1$ & $2$ & $2$ & $4$ & $7$\\\hline
$\left\{  \deg y_{t_{1}},\cdots,\deg y_{t_{k}}\right\}  $ & $\{6\}$ &
$\{6,8\}$ & $\{6,8\}$ & $\{6,8,10,18\}$ & $\{6,8,10,12,18,20,30\}$\\\hline
$\{tor(y_{t_{i}})\}$ & $\{2\}$ & $\{2,3\}$ & $\{2,3\}$ & $\{2,3,2,2\}$ &
$\{2,3,2,5,2,3,2\}$\\\hline
$\{cl(y_{t_{i}})\}$ & $\{2\}$ & $\{2,3\}$ & $\{2,3\}$ & $\{2,3,2,2\}$ &
$\{8,3,4,5,2,3,2\}$\\\hline
\end{tabular}
\\\hline
\end{tabular}
}
\end{center}

By the contents of third rows, the degree sequence $\left\{  \deg y_{t_{1}%
},\cdots,\deg y_{t_{k}}\right\}  $ is strictly increasing for any $G$. For
this reason we shall adopt the convention that $t_{i}:=\deg y_{t_{i}}$, $1\leq
i\leq k$, throughout the remaining parts of this paper.\hfill$\square
$\noindent

\bigskip

Fix, once and for all, for each Lie group $G$ a set of special Schubert
classes $\left\{  y_{t_{1}},\cdots,y_{t_{k}}\right\}  $ on $G/T$, where $\deg
y_{t}=t$. By Lemma 1.4 the inclusions $\omega_{i},y_{t_{i}}\in H^{\ast}(G/T)$
extend to an epimorphism

\begin{quote}
$I:\mathbb{Z}[\omega_{1},\cdots,\omega_{n},y_{t_{1}},\cdots,y_{t_{k}%
}]\rightarrow H^{\ast}(G/T)$.
\end{quote}

\noindent As the ring $\mathbb{Z}[\omega_{1},\cdots,\omega_{n},y_{t_{1}%
},\cdots,y_{t_{k}}]$ is graded by $\deg\omega_{i}=2$ and $\deg y_{t}=t$, we
may use $l(g)$ to denote the degree of a homogeneous polynomial $g$ therein.
According to the Hilbert basis theorem, there exists a minimal set of
homogeneous polynomials $\left\{  R_{1},\ldots,R_{m}\right\}  $ such that
$\ker I$ is the ideal $\left\langle R_{1},\ldots,R_{m}\right\rangle $
generated by $R_{1},\ldots,R_{m}$. In particular, $I$ induces an isomorphism
of rings

\begin{quote}
$H^{\ast}(G/T)=\mathbb{Z}[\omega_{1},\ldots,\omega_{n},y_{t_{1}}%
,\ldots,y_{t_{k}}]/\left\langle R_{1},\ldots,R_{m}\right\rangle $
\end{quote}

\noindent that is compatible with the Borel-Chevalley formula (1.4) of the
real cohomology $H^{\ast}(G/T;\mathbb{R})$. The following result, established
in \cite[Theorem 1.2]{DZ2} (see also \cite[Theorem 5.5]{DZ3}), summaries
useful properties of the polynomials $R_{1},\ldots,R_{m}$.

\bigskip

\noindent\textbf{Theorem A.}\textsl{\ For each Lie group }$G$\textsl{ there
exist homogenous polynomials }$f_{i},g_{i},e_{j}\in\ker I$\textsl{\ with
}$1\leq i\leq k,1\leq j\leq h$\textsl{,} \textsl{such that }

\begin{enumerate}
\item[(1.7)] $H^{\ast}(G/T)=\mathbb{Z}[\omega_{1},\ldots,\omega_{n},y_{t_{1}%
},\ldots,y_{t_{k}}]/\left\langle f_{i},g_{i},e_{j}\right\rangle _{1\leq i\leq
k,1\leq j\leq h}$\textsl{,}
\end{enumerate}

\noindent\textsl{where}

\begin{quote}
\textsl{i) for each }$1\leq i\leq k$ \textsl{the pair\ of polynomials
}$\left\{  f_{i},g_{i}\right\}  $ \textsl{is related to the special Schubert
class }$y_{t_{i}}$\textsl{ in the fashions}

$\quad\quad f_{i}$\textsl{\ }$=$\textsl{\ }$p_{i}\cdot y_{t_{i}}-a_{i}%
$\textsl{, }$g_{i}=y_{t_{i}}^{r_{i}}-b_{i}$\textsl{, }

\noindent\textsl{where }$p_{i}=tor(y_{t_{i}})$\textsl{,} $r_{i}=cl(y_{t_{i}}%
)$\textsl{\ and }$a_{i},b_{i}\in\left\langle \omega_{1},\ldots,\omega
_{n}\right\rangle $\textsl{;}

\textsl{ii) }$e_{j}\in\left\langle \omega_{1},\ldots,\omega_{n}\right\rangle
$\textsl{, }$1\leq j\leq h$\textsl{;}

\textsl{iii) the degree map} $l:\left\{  g_{i},e_{j}\right\}  _{1\leq i\leq
k,1\leq j\leq h}\rightarrow\mathbb{Z}$ \textsl{surjects onto} \textsl{the set}
$2\cdot q(G)$\textsl{, which is also one to one with the only exception}
\end{quote}

\begin{enumerate}
\item[(1.8)] $l^{-1}(2l_{8})=\{g_{4},g_{6},g_{7}\}$\textsl{ for} $G=E_{8}%
$\textsl{.}\hfill$\square$
\end{enumerate}

\bigskip

For a special Schubert class $y_{t_{i}}$ the pair $\left\{  f_{i}%
,g_{i}\right\}  $ of relations is transparent by (1.6). In addition, we note
by the contents in the fourth rows of Table 2 that

\bigskip

\noindent\textbf{Corollary 1.7.} \textsl{The torsion index }$p_{i}$ \textsl{of
a\textsl{ special }Schubert class }$y_{t_{i}}$ \textsl{is a prime, which takes
its value only in} $\{2,3,5\}$\textsl{.}\hfill$\square$

\bigskip

The dimensions of the simply Lie groups $G$ are well-known to be

\begin{center}%
\begin{tabular}
[c]{l|llllllll}\hline\hline
$G$ & $SU(n)$ & $Spin(n)$ & $\qquad Sp(n)$ & $G_{2}$ & $F_{4}$ & $E_{6}$ &
$E_{7}$ & $E_{8}$\\\hline
$\dim G$ & $n^{2}-1$ & $\frac{n(n-1)}{2}$ & $\ n(2n+1)$ & $\ 14$ & $52$ & $78$
& $133$ & $248$\\\hline\hline
\end{tabular}
.
\end{center}

\noindent Property iii) of Theorem A indicates a relationship between $\dim G$
and the degrees of the polynomials $g_{i},e_{j}$ in (1.7). In particular, we have

\begin{quote}
$k+h=n$ for $G\neq E_{8}$, but $k+h=10>8$ for $G=E_{8}$,
\end{quote}

\noindent while the routine formula due to Chevalley \cite{Ch1}

\begin{quote}
$\dim G=(2l_{1}-1)+\cdots+(2l_{n}-1)$
\end{quote}

\noindent goes over to

\bigskip

\noindent\textbf{Corollary 1.8.} $\dim G=\underset{1\leq j\leq h}{\Sigma
}(l(e_{j})-1)+\Sigma(l(g_{i})-1)$\textsl{, where the second sum ranges over
}$1\leq i\leq k$\textsl{, with the only exception that }$i\neq4,7$\textsl{ for
}$G=E_{8}.$\hfill$\square$\noindent

\bigskip

\noindent\textbf{Remark 1.9. }For a Lie group $G\neq E_{8}$ the polynomials
$e_{i},f_{j},g_{j}$ in (1.7) can be shown to be algebraically independent in
the numerator. In contrast, if $G=E_{8}$ one has

\begin{quote}
$\deg g_{4}=\deg g_{6}=\deg g_{7}=60$
\end{quote}

\noindent by (1.8), while there exists a polynomial of the form

\begin{quote}
$\phi=2y_{12}^{5}-y_{20}^{3}+y_{30}^{2}+\beta$ with $\beta\in\left\langle
\omega_{1},\cdots,\omega_{8}\right\rangle $
\end{quote}

\noindent such that the three polynomials $g_{4},g_{6}$ and $g_{7}$ are
related by the formula (see \cite[(6.1), (6.2)]{DZ2})

\begin{enumerate}
\item[(1.9)] $\left\{
\begin{tabular}
[c]{l}%
$g_{4}=-12\phi+5y_{12}^{4}f_{4}-4y_{20}^{2}f_{6}+6y_{30}f_{7}$;\\
$g_{6}=-10\phi+4y_{12}^{4}f_{4}-3y_{20}^{2}f_{6}+5y_{30}f_{7}$;\\
$g_{7}=15\phi-6y_{12}^{4}f_{4}+5y_{20}^{2}f_{6}-7y_{30}f_{7}$.
\end{tabular}
\right.  $
\end{enumerate}

\noindent This phenomenon will cause a few additional concerns for the case
$G=E_{8}$ in our unified approach to $H^{\ast}(G)$.\hfill$\square$

\subsection{The generators of the ring $H^{\ast}(G)$\hfill}

Suppose that $G$ is a Lie group of rank $n$ with a fixed set $\left\{
y_{t_{1}},\cdots,y_{t_{k}}\right\}  $ of special Schubert classes on
$G/T$.\textsl{ }Inputting the formula (1.7) into (1.3) one gets a concise
expression of the bi-graded ring

\begin{quote}
$E_{2}^{\ast,\ast}(G)=\frac{\mathbb{Z}[\omega_{1},\ldots,\omega_{n},y_{t_{1}%
},\cdots,y_{t_{k}}]}{\left\langle f_{i},g_{i},e_{j},\text{ }1\leq i\leq
k;1\leq j\leq h\right\rangle }\otimes\Lambda(\lambda_{1},\cdots,\lambda_{n}).$
\end{quote}

\noindent It enables us to construct a minimal set of generators of the
cohomology $H^{\ast}(G)$, uniformly for all $G$. We begin by introducing two
maps $\left[  \mathcal{D}\right]  $ and $\kappa$, critical for the construction.

Firstly, in term of the monomials basis $\{\omega_{1}^{k_{1}}\cdots\omega
_{n}^{k_{n}}y^{\alpha},k_{1}+\cdots+k_{n}\geq1\}$ of the ideal $\left\langle
\omega_{1},\cdots,\omega_{n}\right\rangle $ consider the linear map of degree
$-1$

\begin{enumerate}
\item[(1.10)] $\mathcal{D}:\left\langle \omega_{1},\cdots,\omega
_{n}\right\rangle \rightarrow E_{2}^{\ast,1}(G)=H^{\ast}(G/T)\otimes
\Lambda^{1}(\lambda_{1},\cdots,\lambda_{n})$
\end{enumerate}

\noindent defined by

\begin{quote}
$\mathcal{D}(\omega_{1}^{k_{1}}\cdots\omega_{n}^{k_{n}}\cdot y^{\alpha
})=I(\omega_{s}^{k_{s}-1}\cdots\omega_{n}^{k_{n}}\cdot y^{\alpha}%
)\otimes\lambda_{s}$,
\end{quote}

\noindent where $s\in\{1,\cdots,n\}$ is the least index such that $k_{s}\geq
1$, $y^{\alpha}$ denotes an arbitrary monomial in $y_{t_{1}},\cdots,y_{t_{k}}%
$. By Lemma 1.1 the map $\mathcal{D}$ fits into the commutative triangle

\begin{quote}
$%
\begin{array}
[c]{ccc}%
\mathcal{D} &  & E_{2}^{\ast,1}(G)=H^{\ast}(G/T)\otimes\Lambda^{1}(\lambda
_{1},\cdots,\lambda_{n})\\
& \nearrow & \downarrow d_{2}\\
\left\langle \omega_{1},\cdots,\omega_{n}\right\rangle  & \overset{I^{\prime}%
}{\rightarrow} & E_{2}^{\ast,0}(G)=H^{\ast}(G/T)\text{,}%
\end{array}
$
\end{quote}

\noindent where $I^{\prime}$ is the restriction of the empimorphism $I$ to
$\left\langle \omega_{1},\cdots,\omega_{n}\right\rangle $. Observe that if
$a\in\left\langle \omega_{1},\cdots,\omega_{n}\right\rangle \cap\ker I$ then

\begin{quote}
$d_{2}(\mathcal{D}(a))=I^{\prime}(a)=I(a)=0$,
\end{quote}

\noindent implying $\mathcal{D}(a)\in\ker d_{2}$. For a $d_{2}$%
-cocycle\textsl{ }$\theta\in\ker d_{2}$ write $[\theta]\in E_{3}^{\ast,\ast
}(G)$\textsl{ }to denote its cohomology class. We have shown

\bigskip

\noindent\textbf{Lemma 1.10.} \textsl{The map }$\left[  \mathcal{D}\right]
:\left\langle \omega_{1},\cdots,\omega_{n}\right\rangle \cap\ker I\rightarrow
$\textsl{ }$E_{3}^{\ast,1}(G)$ \textsl{by }$\left[  \mathcal{D}\right]
(a)=\left[  \mathcal{D}(a)\right]  $ \textsl{is linear, and satisfies
}$\left[  \mathcal{D}\right]  (a)\in E_{3}^{\deg a-2,1}(G).$\hfill$\square
$\noindent

\bigskip

Next, the relations $d_{r}(E_{r}^{\ast,1})=0$ for all $r\geq3$ give rise to a
sequence of quotient maps (for details, see (2.4))

\begin{enumerate}
\item[(1.11)] $\kappa:E_{3}^{2k,1}(G)\twoheadrightarrow\cdots
\twoheadrightarrow E_{\infty}^{2k,1}(G)=\mathcal{F}^{2k}(H^{2k+1}(G))\subset
H^{2k+1}(G)$,
\end{enumerate}

\noindent that interprets elements of $E_{3}^{2k,1}(G)$ directly as cohomology
classes of $G$, where $\mathcal{F}$ is the filtration on $H^{\ast}(G)$ induced
by $\pi$. Granted with the maps $\left[  \mathcal{D}\right]  $ and $\kappa$ we
introduce three types of elements of $H^{\ast}(G)$.

\bigskip

\noindent\textbf{Definition 1.11. }Let $\pi^{\ast}:H^{\ast}(G/T)\rightarrow
H^{\ast}(G)$ be the induced map of $\pi$. For a special Schubert classes
$y_{t_{i}}$ on $G/T$ define the \textsl{Schubert cocycle} $x_{t_{i}}$ on $G$ by

\begin{enumerate}
\item[(1.12)] $x_{t_{i}}:=\pi^{\ast}(y_{t_{i}})\in H^{t_{i}}(G)$, $1\leq i\leq
k$.\hfill$\square$
\end{enumerate}

\bigskip

\noindent\textbf{Definition 1.12. }By ii) of Theorem A the polynomials $e_{j}$ satisfies

\begin{quote}
$e_{j}\in$ $\left\langle \omega_{1},\ldots,\omega_{n}\right\rangle \cap\ker
I$, $1\leq j\leq h$.
\end{quote}

\noindent Applying the composition $\kappa\circ\lbrack\mathcal{D}]$ we obtain
the cohomology classes

\begin{enumerate}
\item[(1.13)] $\rho_{l(e_{j})-1}:=\kappa\circ\lbrack\mathcal{D}](e_{j})\in
H^{l(e_{j})-1}(G),$ $1\leq j\leq h$.
\end{enumerate}

Similarly, by i) of Theorem A, the pair\ of polynomials $f_{i},g_{i}\in\ker I$
yields following the element belonging to $\left\langle \omega_{1}%
,\ldots,\omega_{n}\right\rangle \cap\ker I$

\begin{quote}
$y_{t_{i}}^{r_{i}-1}\cdot f_{i}-p_{i}\cdot g_{i}=y_{t_{i}}^{r_{i}-1}\cdot
a_{i}-p_{i}\cdot b_{i}$.
\end{quote}

\noindent Consequently, we obtain the cohomology class

\begin{enumerate}
\item[(1.14)] $\rho_{l(g_{i})-1}:=\kappa\circ\lbrack\mathcal{D}](y_{i}%
^{r_{i}-1}\cdot a_{i}-p_{i}\cdot b_{i})\in H^{l(g_{i})-1}(G)$,
\end{enumerate}

\noindent where $1\leq i\leq k$, with the exception that $i\neq4,7$ if
$G=E_{8}$ (see (1.8)).

For convenience, we refer the classes $\rho_{l(e_{j})-1}$, $\rho_{l(g_{i})-1}$
constructed in (1.13) and (1.14) as \textsl{the primary classes }of the
cohomology $H^{\ast}(G)$.\hfill$\square$

\bigskip

A special Schubert class $y_{t}$ on $G/T$ is called $p$\textsl{-special }if
$tor(y_{t})=p$. Let $D_{1}(G,p)$ be the degree set of the $p$-special Schubert
classes on $G/T$, to be arranged as a sequence with respect to the order $>$
among integers. As examples, one reads from the last two columns of Table 2 that

\begin{center}%
\begin{tabular}
[c]{l||l|l|l|l}\hline
$p$ & $2$ & $3$ & $5$ & $>5$\\\hline\hline
$D_{1}(E_{7},p)$ & $\left\{  6,10,18\right\}  $ & $\left\{  8\right\}  $ &
$\emptyset$ & $\emptyset$\\\hline
$D_{1}(E_{8},p)$ & $\left\{  6,10,18,30\right\}  $ & $\left\{  8,20\right\}  $
& $\left\{  12\right\}  $ & $\emptyset$\\\hline
\end{tabular}
.
\end{center}

\noindent\textbf{Definition 1.13. }For each $t_{i}\in D_{1}(G,p)$ consider the
polynomial $f_{i}$\textsl{\ }$=$\textsl{\ }$p\cdot y_{t_{i}}-a_{i}$ related to
the Schubert class $y_{t_{i}}$. From $f_{i}\in\ker I$ and $a_{i}%
\in\left\langle \omega_{1},\ldots,\omega_{n}\right\rangle $ by i) of Theorem A
one finds that

\begin{quote}
$f_{i}\equiv a_{i}\operatorname{mod}p\in\left\langle \omega_{1},\ldots
,\omega_{n}\right\rangle \cap\ker I_{p}$, where $I_{p}=I\operatorname{mod}p$.
\end{quote}

\noindent Thus, letting $\kappa^{\prime}$ and $\left[  D\right]  ^{\prime}$ be
respectively the $\mathbb{F}_{p}$-analogues of the maps $\kappa$ and $\left[
\mathcal{D}\right]  $ (see in Section \S 3.2), we obtain the
$\operatorname{mod}p$ cohomology class of $G$

\begin{enumerate}
\item[(1.15)] $\theta_{t_{i}-1}:=\kappa^{\prime}\circ\lbrack\mathcal{D}%
^{\prime}](a_{i})\in$ $H^{t_{i}-1}(G;\mathbb{F}_{p})$.
\end{enumerate}

A subsequence $J\subseteq D_{1}(G,p)$ is called $p$\textsl{-monotone} if its
cardinality $\left\vert J\right\vert \geq2$. For such a sequence $J$ we
formulate, in the sequel to (1.15), the $p$-torsion element

\begin{enumerate}
\item[(1.16)] $\mathcal{C}_{J}:=\beta_{p}(\theta_{I})\in H^{\ast}(G)$ with
$\theta_{J}=\Pi_{s\in J}\theta_{s-1}$,
\end{enumerate}

\noindent where $\beta_{p}:H^{r}(G;\mathbb{F}_{p})\rightarrow H^{r+1}(G)$ is
the Bockstein homomorphism.\hfill$\square$

\bigskip

\noindent\textbf{Example 1.14. }Since\textbf{ }the set $q(G)$ consists
distinct integers by of Table 1, the totality of the primary classes
$\{\rho_{l(e_{j})-1},\rho_{l(g_{i})-1}\}$ of $G$ can be rephrased as
$\{\rho_{2l-1},l\in q(G)\}$ by iii) of Theorem A. In particular, by Corollary 1.8,

\begin{enumerate}
\item[(1.17)] $\dim G=\Sigma_{l\in q(G)}\deg\rho_{2l-1}$.
\end{enumerate}

In addition, in terms of the invariants of the group $G$ given by Table 2, one
can enumerate all the Schubert cocycles $x_{s}$, and the torsion elements
$\mathcal{C}_{J}$ of $H^{\ast}(G)$. For example, if $G=E_{7}$, these classes
are $\left\{  x_{6},x_{8},x_{10},x_{18}\right\}  $ and $\mathcal{C}_{J}$,
where $J\subseteq\left\{  6,10,18\right\}  $ are $2$\textsl{-}monotone.\hfill
$\square$

\subsection{The cohomology of exceptional Lie groups}

The integral cohomology of any finite $CW$--complex $X$ admits the decomposition

\begin{enumerate}
\item[(1.18)] $H^{\ast}(X)=\mathcal{F}(X)\oplus_{p}\tau_{p}(X)$,
\end{enumerate}

\noindent where $\mathcal{F}(X):=H^{\ast}(X)/TorH^{\ast}(X)$ is \textsl{the
free part} of $H^{\ast}(X)$, the direct sum $\oplus$ is over all primes
$p\geq2$, and where $\tau_{p}(X)$ is the $p$-\textsl{primary component }of
$H^{\ast}(X)$ defined by

\begin{quote}
\textsl{\ }$\tau_{p}(X):=\{x\in H^{\ast}(X)\mid p^{r}x=0,$ $r\geq1\}$.
\end{quote}

Given a sequence of graded elements $z_{1},\cdots,z_{m}$ denote by
$\Delta(z_{1},\cdots,z_{m})$ \footnote{This notion $\Delta(z_{1},\cdots
,z_{m})$ is due to Borel, called the group in the simple system of generators
$z_{1},\cdots,z_{m}$.} the free $\mathbb{Z}$-module with the monomial basis
$z_{1}^{\varepsilon_{1}}\cdots z_{m}^{\varepsilon_{m}}$, where $\varepsilon
_{i}\in\{0,1\}$ and $z_{i}^{0}=1$. In addition, if $B=B^{0}\oplus B^{1}\oplus
B^{2}\oplus\cdots$ is a graded group (resp. ring), we use $B^{+}$ to denote
its subgroup (resp. subring) $B^{1}\oplus B^{2}\oplus\cdots$.

The integral cohomologies $H^{\ast}(G)$ have been computed for $G=SU(n)$,
$Sp(n)$ by Borel \cite{B1}, and for $G=Spin(n)$ by Pittie \cite{P}. Our main
result presents the cohomologies of the exceptional Lie groups, merely using
integral classes $x_{t}$, $\varrho_{2l-1}$ and $\mathcal{C}_{J}$ constructed
in Section \S 1.2.

\bigskip

\noindent\textbf{Theorem B.} \textsl{The integral cohomology of the
exceptional Lie groups are:}

\begin{enumerate}
\item[i)] $H^{\ast}(G_{2})=\Delta(\varrho_{3})\otimes\Lambda(\varrho
_{11})\oplus\tau_{2}(G_{2})$\textsl{, where}

$\qquad\tau_{2}(G_{2})=\mathbb{F}_{2}[x_{6}]^{+}/\left\langle x_{6}%
^{2}\right\rangle \otimes\Delta(\varrho_{3})$\textsl{, }

\textsl{and where the generators are subject to the relations}

$\qquad\varrho_{3}^{2}=x_{6}$\textsl{,} $\varrho_{11}\cdot x_{6}=0$\textsl{.}

\item[ii)] $H^{\ast}(F_{4})=\Delta(\varrho_{3})\otimes\Lambda(\varrho
_{11},\varrho_{15},\varrho_{23})\oplus\tau_{2}(F_{4})\oplus\tau_{3}(F_{4}%
)$\textsl{, where}

$\qquad\tau_{2}(F_{4})=\mathbb{F}_{2}[x_{6}]^{+}/\left\langle x_{6}%
^{2}\right\rangle \otimes\Delta(\varrho_{3})\otimes\Lambda(\varrho
_{15},\varrho_{23})$\textsl{,}

$\qquad\tau_{3}(F_{4})=\mathbb{F}_{3}[x_{8}]^{+}/\left\langle x_{8}%
^{3}\right\rangle \otimes\Lambda(\varrho_{3},\varrho_{11},\varrho_{15}%
)$\textsl{,}

\textsl{and} \textsl{where the generators are subject to the relations}

$\qquad\varrho_{3}^{2}=x_{6}$\textsl{, } $\varrho_{11}\cdot x_{6}=0$\textsl{,
}$\varrho_{23}\cdot x_{8}=0$\textsl{.}

\item[iii)] $H^{\ast}(E_{6})=\Delta(\varrho_{3})\otimes\Lambda(\varrho
_{9},\varrho_{11},\varrho_{15},\varrho_{17},\varrho_{23})\oplus\tau_{2}%
(E_{6})\oplus\tau_{3}(E_{6})$\textsl{, where}

$\qquad\tau_{2}(E_{6})=\mathbb{F}_{2}[x_{6}]^{+}/\left\langle x_{6}%
^{2}\right\rangle \otimes\Delta(\varrho_{3})\otimes\Lambda(\varrho_{9}%
,\varrho_{15},\varrho_{17},\varrho_{23})$\textsl{,}

$\qquad\tau_{3}(E_{6})=\mathbb{F}_{3}[x_{8}]^{+}/\left\langle x_{8}%
^{3}\right\rangle \otimes\Lambda(\varrho_{3},\varrho_{9},\varrho_{11}%
,\varrho_{15},\varrho_{17})$\textsl{,}

\textsl{and where the generators are subject to the relations}

$\qquad\varrho_{3}^{2}=x_{6}$\textsl{, }$\varrho_{11}\cdot x_{6}=0$\textsl{,
}$\varrho_{23}\cdot x_{8}=0$\textsl{.}

\item[iv)] $H^{\ast}(E_{7})=\Delta(\varrho_{3})\otimes\Lambda(\varrho
_{11},\varrho_{15},\varrho_{19},\varrho_{23},\varrho_{27},\varrho
_{35})\underset{p=2,3}{\oplus}\tau_{p}(E_{7})$\textsl{, where}

$\qquad\tau_{2}(E_{7})=\frac{\mathbb{F}_{2}[x_{6},x_{10},x_{18},\mathcal{C}%
_{I}]^{+}}{\left\langle x_{6}^{2},x_{10}^{2},x_{18}^{2},\mathcal{R}%
_{J},\mathcal{D}_{K},\mathcal{S}_{H,L}\right\rangle }\otimes\Delta(\varrho
_{3})\otimes\Lambda(\varrho_{15},\varrho_{23},\varrho_{27})$

\textsl{with }$I,J,H,L\subseteq\{6,10,18\}$\textsl{, }$\left\vert I\right\vert
,\left\vert J\right\vert ,\left\vert L\right\vert \geq2$\textsl{, and
}$K=\{6,10,18\}$\textsl{;}

$\qquad\tau_{3}(E_{7})=\frac{\mathbb{F}_{3}[x_{8}]^{+}}{\left\langle x_{8}%
^{3}\right\rangle }\otimes\Lambda(\varrho_{3},\varrho_{11},\varrho
_{15},\varrho_{19},\varrho_{27},\varrho_{35})$\textsl{,}

\textsl{and} \textsl{where the generators are subject to the relations}

$\qquad\varrho_{3}^{2}=x_{6}$\textsl{, }$\varrho_{23}\cdot x_{8}=0,$

\qquad$\mathcal{H}_{i,K}\in\tau_{2}(E_{7})$ \textsl{with} $i\in\{6,10,18\}$%
\textsl{,} $K\subseteq\{6,10,18\}$.

\item[v)] $H^{\ast}(E_{8})=\Delta(\varrho_{3},\varrho_{15},\varrho
_{23})\otimes\Lambda(\varrho_{27},\varrho_{35},\varrho_{39},\varrho
_{47},\varrho_{59})\underset{p=2,3,5}{\oplus}\tau_{p}(E_{8})$\textsl{, where}

$\qquad\tau_{2}(E_{8})=\frac{\mathbb{F}_{2}[x_{6},x_{10},x_{18},x_{30}%
,\mathcal{C}_{I}]^{+}}{\left\langle x_{6}^{8},x_{10}^{4},x_{18}^{2},x_{30}%
^{2},\mathcal{R}_{J},\mathcal{D}_{K}\mathcal{S}_{H,L}\right\rangle }%
\otimes\Delta(\varrho_{3},\varrho_{15},\varrho_{23})\otimes\Lambda
(\varrho_{27})$

\textsl{with} $I$\textsl{,}$J,K,H,L\subseteq\{6,10,18,30\}$\textsl{ and
}$\left\vert I\right\vert ,\left\vert J\right\vert ,\left\vert L\right\vert
\geq2$\textsl{, }$\left\vert K\right\vert \geq3$\textsl{,}

$\qquad\tau_{3}(E_{8})=\frac{\mathbb{F}_{3}[x_{8},x_{20},\mathcal{C}%
_{\{8,20\}}]^{+}}{\left\langle x_{8}^{3},x_{20}^{3},x_{8}^{2}x_{20}%
^{2}\mathcal{C}_{\{8,20\}},\mathcal{C}_{\{8,20\}}^{2}\right\rangle }%
\otimes\Lambda(\varrho_{3},\varrho_{15},\varrho_{27},\varrho_{35},\varrho
_{39},\varrho_{47})$\textsl{,}

$\qquad\tau_{5}(E_{8})=\frac{\mathbb{F}_{5}[x_{12}]^{+}}{\left\langle
x_{12}^{5}\right\rangle }\otimes\Lambda(\varrho_{3},\varrho_{15},\varrho
_{23},\varrho_{27},\varrho_{35},\varrho_{39},\varrho_{47})$\textsl{,}

\textsl{and where the generators are subject to the relations}

$\qquad\varrho_{3}^{2}=x_{6}$\textsl{, }$\varrho_{15}^{2}=x_{30}$\textsl{,
}$\varrho_{23}^{2}=x_{6}^{6}x_{10}$\textsl{, }$\varrho_{59}\cdot x_{12}=0$;

$\qquad\varrho_{23}\cdot x_{8}=0$, $\rho_{23}\cdot x_{20}=x_{8}^{2}%
\cdot\mathcal{C}_{\left\{  8,20\right\}  }$, $\rho_{23}\cdot\mathcal{C}%
_{\left\{  8,20\right\}  }=0$;

$\qquad\rho_{59}\cdot x_{20}=0$, $\rho_{59}\cdot x_{8}=x_{20}^{2}%
\mathcal{C}_{\left\{  8,20\right\}  }$, $\rho_{59}\cdot\mathcal{C}_{\left\{
8,20\right\}  }=0$;

\qquad$\mathcal{H}_{i,K}\in\tau_{2}(E_{8})$ \textsl{with} $i\in\{6,10,18,30\}$%
\textsl{,} $K\subseteq\{6,10,18,30\}$.
\end{enumerate}

\noindent\textsl{In addition, the relations of the types} $\mathcal{R}%
_{J},\mathcal{D}_{K},\mathcal{S}_{H,L}$ \textsl{and} $\mathcal{H}_{i,K}%
$\textsl{,} \textsl{concerning only} \textsl{with} \textsl{the ideals }%
$\tau_{2}(E_{7})$ \textsl{and} $\tau_{2}(E_{8})$\textsl{, are too many to be
listed explicitly. Instead, their general formulae are stated in (4.4), (4.5)
and (4.7), respectively.}

\bigskip

The remaining sections of the paper are arranged as following. Section \S 2
develops properties of the cohomologies of certain Koszul complexes, useful to
solve the extension problem from $E_{3}^{\ast,\ast}(G)$ to $H^{\ast}(G)$.
Section \S 3 recalls from \cite{DZ1} a presentation of the algebra $H^{\ast
}(G;\mathbb{F}_{p})$, and calculates the $\operatorname{mod}p$ Bockstein
cohomology $H_{\beta}^{\ast}(G;\mathbb{F}_{p})$. In Section \S 4 we determine
the structures of three basic components of the integral cohomology $H^{\ast
}(G)$: the subring $\operatorname{Im}\pi^{\ast}$, the free part $\mathcal{F}%
(G)$, and the torsion ideals $\tau_{p}(G)$. As applications, Theorem B is
proved in Section \S 4.4.

The study of the cohomology of Lie groups has a long and outstanding history.
Section \S 5 recalls early works on the subject, make comparison of our
approach with the classical ones. Along the way we illustrate the perspective
and necessity of the conceptual developments in the present work.

\section{The cohomology of certain Koszul complexes}

Let $A$ be a graded ring (or algebra), and let $\{z_{1},\ldots,z_{k}\}$ be a
sequence of homogeneous elements of $A$. The \textsl{Koszul complex}
associated to the pair $\{A;(z_{1},\ldots,z_{k})\}$, written $K(A;z_{1}%
,\ldots,z_{k})$, is the cochain complex $\left\{  C,\delta\right\}  $, where

\begin{quote}
i) $C$ is the free $A$-module $A\otimes\Delta(\theta_{1},\cdots,\theta_{k})$
graded $\deg\theta_{t}=\deg z_{t}-1$;

ii) $\delta$ is the antiderivation of degree $1$ on $C$ defined by

$\qquad\delta(1\otimes\theta_{t})=z_{t}\otimes1$ and $\delta(z\otimes1)=0$,
$z\in A$.
\end{quote}

\noindent With $\delta\circ\delta=0$ the cohomology $H^{\ast}(K(A;z_{1}%
,\ldots,z_{k}))$ of the complex $\left\{  C,\delta\right\}  $ is the graded
quotient group $\ker\delta/\operatorname{Im}\delta$.

By Lemma 1.1, the second page $\left\{  E_{2}^{\ast,\ast}(G),d_{2}\right\}  $
of the Serre spectral sequence of $\pi$ is the Koszul complex $K(H^{\ast
}(G/T);\omega_{1},\ldots,\omega_{n})$ whose cohomology is the third page
$E_{3}^{\ast,\ast}(G)$. In this section we relate the groups $E_{3}^{\ast
,0}(G)$ and $E_{3}^{\ast,1}(G)$ with appropriate subgroups of the cohomology
$H^{\ast}(G)$, which helps to solve the extension problem from $E_{3}%
^{\ast,\ast}(G)$ to $H^{\ast}(G)$. We introduce also the notion of the Koszul
complex associated to a polynomial algebra, useful to formulate of the ideal
$\tau_{p}(G)$ in Section \S 4.

\subsection{The Koszul complex $\left\{  E_{2}^{\ast,\ast}(G),d_{2}\right\}
$}

Let $\{E_{r}^{\ast,\ast}(G),d_{r}\}$ be the Serre spectral sequence of the
torus fibration $\pi$, and let $\mathcal{F}^{p}$ be the filtration on
$H^{\ast}(G)$ defined by $\pi$. That is (\cite[P.146]{Mc})

\begin{center}
$0=\mathcal{F}^{r+1}(H^{r}(G))\subseteq\mathcal{F}^{r}(H^{r}(G))\subseteq
\cdots\subseteq\mathcal{F}^{0}(H^{r}(G))=H^{r}(G)$
\end{center}

\noindent with

\begin{enumerate}
\item[(2.1)] $E_{\infty}^{p,q}(G)=\mathcal{F}^{p}(H^{p+q}(G))/\mathcal{F}%
^{p+1}(H^{p+q}(G))$.
\end{enumerate}

The routine property $d_{r}(E_{r}^{\ast,0}(G))=0$, $r\geq2$, gives rise to a
sequence of quotient maps

\begin{center}
$H^{r}(G/T)=E_{2}^{r,0}\twoheadrightarrow E_{3}^{r,0}\twoheadrightarrow
\cdots\twoheadrightarrow E_{\infty}^{r,0}=\mathcal{F}^{r}(H^{r}(G))\subseteq
H^{r}(G)$
\end{center}

\noindent whose composition agrees with $\pi^{\ast}:H^{\ast}(G/T)\rightarrow
H^{\ast}(G)$ \cite[P.147]{Mc}. We may therefore reserve $\pi^{\ast}$ for the epimorphism:

\begin{enumerate}
\item[(2.2)] $\pi^{\ast}:E_{3}^{\ast,0}(G)\twoheadrightarrow\mathcal{F}%
^{r}(H^{r}(G))\subset H^{\ast}(G)$.
\end{enumerate}

\noindent In addition, since for any special Schubert class $y_{t_{i}}$ the
class $y_{t_{i}}\otimes1\in E_{2}^{\ast,0}(G)$ is $d_{2}$-closed by Lemma 1.1,
we may reserve $y_{t_{i}}$ to simplify $[y_{t_{i}}\otimes1]\in E_{3}^{\ast
,0}(G)$.

\bigskip

\noindent\textbf{Lemma 2.1.} \textsl{If }$\left\{  y_{t_{1}},\cdots,y_{t_{k}%
}\right\}  $\textsl{ is a set of special Schubert classes on }$G/T$\textsl{,
then}

\begin{quote}
$E_{3}^{\ast,0}(G)=\mathbb{Z}[y_{t_{1}},\cdots,y_{t_{k}}]/\left\langle
p_{i}y_{t_{i}},y_{t_{i}}^{r_{i}}\right\rangle _{1\leq i\leq k}$,
\end{quote}

\noindent\textsl{where} $p_{i}=tor(y_{t_{i}})$\textsl{,} $r_{i}=cl(y_{t_{i}})$.

\bigskip

\noindent\textbf{Proof. }The group $E_{3}^{\ast,0}(G)$ is the cokernel of the differential

\begin{quote}
$d_{2}:E_{2}^{\ast,1}(G)=H^{\ast}(G/T)\otimes\Lambda^{1}(t_{1},\cdots
,t_{n})\rightarrow E_{2}^{\ast,0}(G)=H^{\ast}(G/T)$.
\end{quote}

\noindent Since $\operatorname{Im}d_{2}=\left\langle \omega_{1},\cdots
,\omega_{n}\right\rangle _{c}\subseteq H^{\ast}(G/T)$ by Lemma 1.1, we obtain
the formula of $E_{3}^{\ast,0}(G)$ by setting $\omega_{1}=\cdots=\omega_{n}=0$
in (1.7).\hfill$\square$

\bigskip

As the group $E_{r}^{p,q}(G)$ is a subquotient of $E_{r-1}^{p,q}(G)$, the fact
$H^{2s+1}(G/T)=0$ by Lemma 1.3 implies that

\begin{enumerate}
\item[(2.3)] $E_{r}^{2s+1,q}(G)=0$ for all $s,q\geq0$ and $r\geq2$.
\end{enumerate}

\noindent In particular, taking $(r,q)=(\infty,0)$ one gets by (2.1) that

\begin{quote}
$\mathcal{F}^{2s+1}(H^{2s+1}(G))=\mathcal{F}^{2s+2}(H^{2s+1}(G))=0$.
\end{quote}

\noindent It implies, again by (2.1), that

\begin{quote}
$E_{\infty}^{2s,1}(G)=\mathcal{F}^{2s}(H^{2s+1}(G))\subset H^{2s+1}(G)$.
\end{quote}

\noindent Combining this with $d_{r}(E_{r}^{\ast,1})=0$ for $r\geq3$ gives
rise to the sequence of empimorphisms that interprets elements of
$E_{3}^{2k,1}(G)$ directly as cohomology classes of $G$ (e.g. (1.11))

\begin{enumerate}
\item[(2.4)] $\kappa:E_{3}^{\ast,1}(G)\twoheadrightarrow E_{4}^{\ast
,1}(G)\twoheadrightarrow\cdots\twoheadrightarrow E_{\infty}^{\ast
,1}(G)=\mathcal{F}^{2s}(H^{2s+1}(G))\subset H^{\ast}(G)$.
\end{enumerate}

With respect to the product inherited from that on $E_{2}^{\ast,\ast}(G)$ the
third page $E_{3}^{\ast,\ast}(G)$ is a bi-graded ring \cite[P.668]{Wh}.
Concerning the relationship between the two maps $\pi^{\ast}$ and $\kappa$ we
prove the next two results.

\bigskip

\noindent\textbf{Lemma 2.2. }\textsl{For any }$\rho\in E_{3}^{2s,1}%
(G)$\textsl{\ one has} $\kappa(\rho)^{2}\in\operatorname{Im}\pi^{\ast}\cap
\tau_{2}(G)$\textsl{.}

\bigskip

\noindent\textbf{Proof.} For an $\rho\in E_{3}^{2s,1}(G)$ the obvious relation
$\rho^{2}=0$ in

\begin{quote}
$E_{3}^{4s,2}(G)\twoheadrightarrow E_{\infty}^{4s,2}(G)=\mathcal{F}%
^{4s}(H^{4s+2}(G))/\mathcal{F}^{4s+1}(H^{4s+2}(G))$
\end{quote}

\noindent implies that $\kappa(\rho)^{2}\in\mathcal{F}^{4s+1}(H^{4s+2}(G))$. With

\begin{quote}
$\mathcal{F}^{4s+1}(H^{4s+2}(G))/\mathcal{F}^{4s+2}(H^{4s+2}(G))=E_{\infty
}^{4s+1,1}(G)=0$
\end{quote}

\noindent by (2.3) one finds further by (2.2) that

\begin{quote}
$\kappa(\rho)^{2}\in\mathcal{F}^{4s+2}(H^{4s+2}(G))\subset\operatorname{Im}%
\pi^{\ast}$.
\end{quote}

\noindent On the other hand, since $\deg\kappa(\rho)=2s+1$ we must have
$2\kappa(\rho)^{2}=0$. The proof is completed by $\kappa(\rho)^{2}\in\tau
_{2}(G)$.\hfill$\square$

\bigskip

For a prime $p$ consider the Bockstein exact sequence associated to the short
exact sequence $0\rightarrow\mathbb{Z}\rightarrow\mathbb{Z}\rightarrow
\mathbb{F}_{p}\rightarrow0$ of coefficients groups

\begin{enumerate}
\item[(2.5)] $\cdots\rightarrow H^{r}(G)\overset{\cdot p}{\rightarrow}%
H^{r}(G)\overset{r_{p}}{\rightarrow}H^{r}(G;\mathbb{F}_{p})\overset{\beta_{p}%
}{\rightarrow}H^{r+1}(G)\rightarrow\cdots$,
\end{enumerate}

\noindent where $\beta_{p}$ is the Bockstein homomorphism, $r_{p}$ is the
$\operatorname{mod}p$ reduction. Since $H^{\ast}(G)$ (resp. $H^{\ast
}(G;\mathbb{F}_{p})$) can be considered as a module over its subring
$\operatorname{Im}\pi^{\ast}$, we may regard (2.5) as an exact sequence in the
$\operatorname{Im}\pi^{\ast}$-modules. On the other hand, since the cohomology
$H^{\ast}(G/T)$ is torsion free, we have the short exact sequence of $d_{2}$-complexes

\begin{enumerate}
\item[(2.6)] $0\rightarrow E_{2}^{\ast,\ast}(G)\overset{\cdot p}{\rightarrow
}E_{2}^{\ast,\ast}(G)\overset{r_{p}}{\rightarrow}E_{2}^{\ast,\ast
}(G;\mathbb{F}_{p})\rightarrow0$,
\end{enumerate}

\noindent whose cohomological exact sequence contains the section

\begin{enumerate}
\item[(2.7)] $\cdots\rightarrow E_{3}^{\ast,1}(G)\overset{\cdot p}%
{\rightarrow}E_{3}^{\ast,1}(G)\overset{r_{p}}{\rightarrow}E_{3}^{\ast
,1}(G;\mathbb{F}_{p})\overset{\overline{\beta}_{p}}{\rightarrow}E_{3}^{\ast
,0}(G)\overset{\cdot p}{\rightarrow}E_{3}^{\ast,0}(G)\rightarrow\cdots$
\end{enumerate}

\noindent that can be viewed as an exact sequence in the $E_{3}^{\ast,0}%
(G)$-modules. Since both $\operatorname{Im}\pi^{\ast}$ and $\operatorname{Im}%
\kappa$ have been identified with appropriate subgroups of $H^{\ast}(G)$ by
(2.2) and (2.4), we have by (2.5) and (2.7) the following exact ladder

\begin{center}%
\begin{tabular}
[c]{llllllll}%
$E_{3}^{s,1}(G)$ & $\overset{\cdot p}{\rightarrow}$ & $E_{3}^{s,1}(G)$ &
$\overset{r_{p}}{\rightarrow}$ & $E_{3}^{s,1}(G;\mathbb{F}_{p})$ &
$\overset{\overline{\beta}_{p}}{\rightarrow}$ & $E_{3}^{s,0}(G)$ &
$\rightarrow\cdots$\\
$\kappa\downarrow$ &  & $\kappa\downarrow$ &  & $\kappa^{\prime}\downarrow$ &
& $\pi^{\ast}\downarrow$ & \\
$H^{s+1}(G)$ & $\overset{\cdot p}{\rightarrow}$ & $H^{s+1}(G)$ &
$\overset{r_{p}}{\rightarrow}$ & $H^{s+1}(G;\mathbb{F}_{p})$ & $\overset
{\beta_{p}}{\rightarrow}$ & $H^{s}(G)$ & $\rightarrow\cdots$%
\end{tabular}
,
\end{center}

\noindent where $\kappa^{\prime}$ is the analogues of $\kappa$ in the
$\operatorname{mod}p$ cohomologies. This shows that

\bigskip

\noindent\textbf{Lemma 2.3. }\textsl{For any }$x\in E_{3}^{\ast,0}(G)$\textsl{
and} $\rho\in E_{3}^{\ast,1}(G)$\textsl{ one has }

\begin{enumerate}
\item[(2.8)] $\kappa(x\cdot\rho)=\pi^{\ast}x\cdot\kappa(\rho)$;

\item[(2.9)] $\beta_{p}\circ\kappa^{\prime}=\pi^{\ast}\circ\overline{\beta
}_{p}$.\hfill$\square$
\end{enumerate}

\subsection{The Koszul complex of a polynomial algebra}

Given a sequence of $r$ evenly graded elements $\left\{  y_{1},\cdots
,y_{r}\right\}  $, and a sequence of $r$ positive integers $\left\{
k_{1},\cdots,k_{r}\right\}  $, consider the (truncated) polynomial algebra
over the field $\mathbb{F}_{p}$

\begin{enumerate}
\item[(2.10)] $A=\mathbb{F}_{p}[y_{1},\cdots,y_{r}]/\left\langle y_{1}^{k_{1}%
},\cdots,y_{r}^{k_{r}}\right\rangle $.
\end{enumerate}

\noindent The Koszul complex\textsl{ }$K(A;y_{1},\ldots,y_{r})$, abbreviated
by $K(A)=\{C,\delta\}$, will be called \textsl{the Koszul complex associated
to} \textsl{the polynomial algebra} $A$. That is

\begin{quote}
$C=A\otimes\Delta(\theta_{1},\cdots,\theta_{k})$, $\delta(1\otimes\theta
_{t})=y_{t}\otimes1$.
\end{quote}

\noindent The cohomology $H^{\ast}(K(A))$ of $K(A)$, as well as the subgroup
$\operatorname{Im}\delta\subset C$, can be explicitly presented.

If $B$ is a graded ring and $\{z_{1},\cdots,z_{k}\}$ is a set of graded
elements, we write $B\cdot\{z_{1},\cdots,z_{k}\}$ to denote the free $B$
module with basis $\{z_{1},\cdots,z_{k}\}$. With this notation the cochain
group $C$ can be rephrased as

\begin{quote}
$C=A\cdot\left\{  1,\theta_{I}\right\}  $, where $I\subseteq\{1,\cdots,r\}$,
$\theta_{I}:=\underset{t\in I}{\Pi}\theta_{t}$.
\end{quote}

\noindent Accordingly, for a multi-index $I=\left\{  i_{1},\cdots
,i_{k}\right\}  \subseteq\{1,\cdots,r\}$ introduce in $C$ the following elements

\begin{enumerate}
\item[(2.11)] $g_{I}:=(\underset{i_{s}\in I}{\Pi}y_{i_{s}}^{k_{i}-1}%
)\cdot\theta_{I}$, $C_{I}=\underset{i_{s}\in I}{\Sigma}(-1)^{s-1}y_{i_{s}%
}\cdot\theta_{I_{s}}$,

$R_{I}:=(\underset{i_{s}\in I}{\Pi}y_{i}^{k_{i}-1})\cdot C_{I}$,
$D_{I}:=\underset{i_{s}\in I}{\Sigma}(-1)^{s-1}y_{i_{s}}\cdot C_{I_{s}}$,
\end{enumerate}

\noindent where $I_{s}=\{i_{1},\cdots,\widehat{i_{s}},\cdots,i_{k}\}$ in
which$\quad\widehat{}$ \quad means omission.

\bigskip\ 

\noindent\textbf{Lemma 2.4.} \textsl{As subspaces} \textsl{of }$C$\textsl{,
the cohomology }$H^{\ast}(K(A))$\textsl{ and the subgroup }$\operatorname{Im}%
\delta$ \textsl{have the presentations}

\begin{enumerate}
\item[(2.12)] $H^{\ast}(K(A))=\Delta(g_{1},\cdots,g_{r})$ \textsl{with}
$g_{i}=y_{i}^{k_{i}-1}\cdot\theta_{i}$\textsl{.}

\item[(2.13)] $\operatorname{Im}\delta=A\cdot\{C_{I}\}\operatorname{mod}%
\{R_{J}$\textsl{, }$D_{K}\},$\textsl{ }
\end{enumerate}

\noindent\textsl{where} $I,J,K\subseteq\{1,\cdots,r\}$\textsl{ with
}$\left\vert I\right\vert ,\left\vert J\right\vert \geq2$\textsl{, and
}$\left\vert K\right\vert \geq3$\textsl{.}

\bigskip

\noindent\textbf{Proof.} The cochain group $C$ has the decomposition
$C=\otimes_{1\leq t\leq r}C_{t}$ in which each factor $C_{t}:=(\mathbb{F}%
_{p}[y_{t}]/\left\langle y_{t}^{k_{t}}\right\rangle )\otimes\Delta(\theta
_{t})$ is an invariant subspace of $\delta$. Thus, letting $\delta_{t}$ be the
restriction of $\delta$ on $C_{t}$, we get by the K\"{u}nneth formula that

\begin{quote}
$H^{\ast}(K(A))=\otimes_{1\leq t\leq r}H^{\ast}(C_{t},\delta_{t})$.
\end{quote}

\noindent We get formula (2.12) by the obvious fact that $H^{\ast}%
(C_{t},\delta_{t})=\Delta(g_{t})$, i.e. the $\mathbb{F}_{p}$-space with basis
$\left\{  1,g_{t}\right\}  $.

To see (2.13) we observe in view of the short exact sequence $0\rightarrow
\ker\delta\rightarrow C\overset{\delta}{\rightarrow}\operatorname{Im}%
\delta\rightarrow0$ that $\delta$ induces a degree $1$ isomorphism in the $A$-modules

\begin{quote}
$\overline{\delta}:C/\ker\delta\overset{\cong}{\rightarrow}\operatorname{Im}%
\delta$.
\end{quote}

\noindent Since, as a vector space over $\mathbb{F}_{p}$, $\ker\delta=H^{\ast
}(K(A))\oplus\operatorname{Im}\delta$ in which

a) $H^{\ast}(K(A))$ has the basis\textsl{ }$\{1,g_{J}\}$ by (2.12);

b) $\operatorname{Im}\delta$ is spanned over $A$ by the elements $C_{K}%
=\delta(\theta_{K})$,

\noindent one gets the presentation

\begin{quote}
$C/\ker\delta\equiv A\cdot\{1,\theta_{I}\}\operatorname{mod}\{1,g_{J},C_{K}\}$,
\end{quote}

\noindent where $I,J,K\subseteq\{1,\cdots,r\}$. Applying the isomorphism
$\overline{\delta}$ to both sides one obtains formula (2.13) from the obvious relations

\begin{quote}
$\overline{\delta}(\theta_{I})=C_{I}$, $\overline{\delta}(g_{J})=R_{J}$,
$\overline{\delta}(C_{K})=D_{K}$,
\end{quote}

\noindent together with the facts that, if either $I$ or $J$ is a singleton
$\{t\}$, then

\begin{quote}
$C_{I}=y_{t}\in A$,$\quad R_{J}=y_{t}^{k_{t}}=0$,
\end{quote}

\noindent and that, if $K=\{k_{1},k_{2}\}$ is a couple, then $D_{K}=y_{k_{1}%
}y_{k_{2}}-y_{k_{1}}y_{k_{2}}=0$.\hfill$\square$

\bigskip

We emphasize at this stage that, if $b=(b_{1},\cdots,b_{r})\subset A$ is a
sequence of homogeneous elements that satisfies the degree constraints\textsl{
}$\deg b_{i}=2\deg\theta_{i}$, $1\leq i\leq r$, then the cochain group
$C=A\otimes\Delta(\theta_{1},\cdots,\theta_{r})$ can be furnished with the
product $\theta_{i}^{2}=b_{i}$ making it the graded polynomial algebra

\begin{quote}
$C=\frac{\mathbb{F}_{p}[y_{1},\cdots,y_{r},\theta_{1},\cdots,\theta_{r}%
]}{\left\langle y_{1}^{k_{1}},\cdots,y_{r}^{k_{r}},\text{ }\theta_{1}%
^{2}-b_{1},\cdots,\theta_{r}^{2}-b_{r}\right\rangle }$,
\end{quote}

\noindent while the complex $K(A)=\{C,\delta\}$\textsl{ }becomes\textsl{ a
graded differential algebra. }Regarding this fact we call the sequence $b$ as
\textsl{an algebra structure} on $K(A)$.

\bigskip

\noindent\textbf{Theorem 2.5}\textsl{. If the Koszul complex }$K(A)$\textsl{
is furnished with an algebra structure }$b$\textsl{, then }$\operatorname{Im}%
\delta$\textsl{ is isomorphic to the (truncated) polynomial algebra}

\begin{enumerate}
\item[(2.14)] $\operatorname{Im}\delta=\frac{\mathbb{F}_{p}[y_{1},\cdots
,y_{r},C_{I}]^{+}}{\left\langle y_{1}^{k_{1}},\cdots,y_{r}^{k_{r}},R_{J}%
,D_{K},S_{H,L}\right\rangle }$\textsl{, }
\end{enumerate}

\noindent\textsl{where} $I,J,K,H,L\subseteq\{1,\cdots,r\}$\textsl{ with
}$\left\vert I\right\vert ,\left\vert J\right\vert ,\left\vert H\right\vert
,\left\vert L\right\vert \geq2$\textsl{, }$\left\vert K\right\vert \geq
3$\textsl{,} \textsl{and where if }$H=\{h_{1},\cdots,h_{k}\}$\textsl{,}

\begin{enumerate}
\item[(2.15)] $S_{H,L}=C_{H}C_{L}+\underset{h_{s}\in H}{\Sigma}(-1)^{s}%
y_{h_{s}}\cdot b_{H_{s}\cap L}\cdot C_{\left\langle H_{s},L\right\rangle }$
\end{enumerate}

\noindent\textsl{in which} $b_{M}=\underset{s\in M}{\Pi}b_{s}\in A$\textsl{,
}$H_{s}=$ $\{h_{1},\cdots,\widehat{h_{s}},\cdots,h_{k}\}$\textsl{ and
}$\left\langle K,L\right\rangle =\{t\in K\cup L\mid t\notin K\cap
L\}$\textsl{.}

\bigskip

\noindent\textbf{Proof. }Granted with the\textbf{ }algebra structure $b$ on
$C$\textsl{ }the product of\textsl{ }any\textbf{ }two\textbf{ }$C_{H}$,
$C_{L}\in C$ can be expanded as an $A$-linear combination in the $C_{I}$'s. Precisely,

\begin{quote}
$C_{H}\cdot C_{L}=\delta(\theta_{H})\cdot\delta(\theta_{L})=\delta
(\delta(\theta_{H})\cdot\theta_{L})$ (by $\delta^{2}=0$)

$=\delta((\Sigma_{h_{s}\in H}(-1)^{s-1}y_{h_{s}}\theta_{H_{s}})\cdot\theta
_{L})$ (since $\theta$ is an antiderivation)

$=\delta(\Sigma_{h_{s}\in H}(-1)^{s-1}y_{j_{s}}(b_{H_{s}\cap L}\cdot
\theta_{\left\langle H_{s},L\right\rangle }))$ (by $\theta_{i}^{2}=b_{i}$)

$=\Sigma_{h_{s}\in H}(-1)^{s-1}y_{h_{s}}b_{H_{s}\cap L}\cdot C_{\left\langle
H_{s},L\right\rangle }$ (by $\delta(\theta_{I})=C_{I}$, $\delta(A)=0$).
\end{quote}

\noindent It follows that, if we introduce in the polynomial algebra
$\mathbb{F}_{p}[y_{1},\cdots,y_{r},C_{I}]^{+}$ the elements $S_{H,L}$ defined
by (2.15), then additively

\begin{quote}
$\frac{\mathbb{F}_{p}[y_{1},\cdots,y_{r},C_{I}]^{+}}{\left\langle y_{1}%
^{k_{1}},\cdots,y_{r}^{k_{r}},R_{J},D_{K}\right\rangle }\operatorname{mod}%
S_{H,L}\equiv A\cdot\{C_{I}\}\operatorname{mod}\{R_{J},D_{K}\}$.
\end{quote}

\noindent That is, we have derived (2.14) from (2.13).\hfill$\square$

\bigskip

\noindent\textbf{Remark 2.6.} If the complex $C$ happens to be a subgroup of
the $\operatorname{mod}p$ cohomology $H^{\ast}(X;\mathbb{F}_{p})$ of a
topological space $X$, and if $p\equiv1\operatorname{mod}2$, then with respect
to the cup product on $H^{\ast}(X;\mathbb{F}_{p})$ we must have $b_{i}%
=\theta_{i}^{2}=0$ for the degree reason $\deg\theta_{i}\equiv
1\operatorname{mod}2$. Consequently, as an algebra over $\mathbb{F}_{p}$,
$C=A\otimes\Lambda(\theta_{1},\cdots,\theta_{r})$, while the formula (2.14) of
$\operatorname{Im}\delta$ becomes

\begin{quote}
$\operatorname{Im}\delta=\frac{\mathbb{F}_{p}[y_{1},\cdots,y_{r},C_{I}]^{+}%
}{\left\langle y_{1}^{k_{1}},\cdots,y_{r}^{k_{r}},R_{J},D_{K},S_{H,L}%
\right\rangle }$ with\textsl{ }$\left\vert I\right\vert ,\left\vert
J\right\vert ,\left\vert H\right\vert ,\left\vert L\right\vert \geq2$\textsl{,
}$\left\vert K\right\vert \geq3$,
\end{quote}

\noindent where $S_{H,L}=C_{H}C_{L}+\underset{H_{s}\cap L=\emptyset}{\Sigma
}(-1)^{s}y_{h_{s}}\cdot C_{H_{s}\cup L}$.\hfill$\square$

\section{The algebra $H^{\ast}(G;\mathbb{F}_{p})$}

Granted with Theorem A we have derived in \cite{DZ1} a presentation of the
$\operatorname{mod}p$ cohomology $H^{\ast}(G/T;\mathbb{F}_{p})$, and
determined accordingly the structure of the algebra $H^{\ast}(G;\mathbb{F}%
_{p})$ as an Hopf algebra over the Steenrod algebra $\mathcal{A}_{p}$. To
facilitate calculation with the integral cohomology $H^{\ast}(G)$ in the
present work, we recall in this section relevant results from \cite{DZ1}.
These will be applied in Theorem 3.7 to obtain closed formulae of the
subgroups $H_{\beta}^{\ast}(G;\mathbb{F}_{p})$, $\operatorname{Im}\delta
_{p}\subset H^{\ast}(G;\mathbb{F}_{p})$, where $H_{\beta}^{\ast}%
(G;\mathbb{F}_{p})$ is the $\operatorname{mod}p$ Bockstein cohomology of $G$,
and $\delta_{p}$ is the Bockstein operator on $H^{\ast}(G;\mathbb{F}_{p})$.

\subsection{The algebra $H^{\ast}(G/T;\mathbb{F}_{p})$}

Suppose that $G$ is a Lie group of rank $n$, and that $p$ is a prime. Recall
from Definition 1.13 that the degree set of the $p$-special Schubert classes
on $G/T$ is

\begin{quote}
$D_{1}(G,p):=\{t_{i}\in\{t_{1},\cdots,t_{k}\},tor(y_{t_{i}})=p\}$.
\end{quote}

\noindent For $G\neq E_{8}$ we let $\overline{D}_{1}(G,p)$ be the complement
of $D_{1}(G,p)$ in $\{t_{1},\cdots,t_{k}\}$, but define $\overline{D}%
_{1}(E_{8},p)$ to be the subsets of the degrees $\{t_{1},\cdots,t_{7}\}$ of
the special Schubert classes on $E_{8}/T$ (see in Table 2)):

\begin{center}%
\begin{tabular}
[c]{l||l|l|l|l}\hline
$p$ & $2$ & $3$ & $5$ & $>5$\\\hline
$\overline{D}_{1}(E_{8},p)$ & $\{8\}$ & $\{6,10,18\}$ & $\{6,8,10,18\}$ &
$\{6,8,10,18,20\}$\\\hline
\end{tabular}
.
\end{center}

Since the ring $H^{\ast}(G/T)$ is torsion free \cite{BS}, one can deduce from
Theorem A a formula of the algebra $H^{\ast}(G/T;\mathbb{F}_{p})$ using the isomorphism

\begin{quote}
$H^{\ast}(G/T;\mathbb{F}_{p})=H^{\ast}(G/T)\otimes\mathbb{F}_{p}$.
\end{quote}

\noindent Precisely, the epimorphism by (1.7)

\begin{quote}
$I_{p}\equiv I\operatorname{mod}p:\mathbb{F}_{p}[\omega_{1},\cdots,\omega
_{n},y_{t_{1}},\cdots,y_{t_{k}}]\rightarrow H^{\ast}(G/T;\mathbb{F}_{p})$.
\end{quote}

\noindent satisfies, by property i) of Theorem A, that

\begin{enumerate}
\item[(3.1)] $I_{p}(f_{i})\equiv-I_{p}(a_{i})$ if $t_{i}\in D_{1}(G,p)$,

\item[(3.2)] $I_{p}(y_{t_{i}})\equiv q_{i}I_{p}(a_{i})$ if $t_{i}\notin
D_{1}(G,p)$,
\end{enumerate}

\noindent where $q_{i}>0$ is the least prime such that $q_{i}\cdot p_{i}%
\equiv1\operatorname{mod}p$. It follows from (1.7) and (3.2) that the map
$I_{p}$ restricts to an epimorphism

\begin{quote}
$\overline{I}_{p}:\mathbb{F}_{p}[\omega_{1},\cdots,\omega_{n},y_{t_{i}%
}]_{t_{i}\in D_{1}(G,p)}\rightarrow H^{\ast}(G/T;\mathbb{F}_{p})$.
\end{quote}

\noindent that satisfies

\begin{quote}
$\ker\overline{I}_{p}=\left\langle a_{i}^{(p)},y_{t_{i}}^{r_{i}}-b_{i}%
^{(p)},g_{s}^{(p)},e_{j}^{(p)}\right\rangle _{t_{i}\in D_{1}(G,p),t_{s}%
\in\overline{D}_{1}(G,p),1\leq j\leq h}$,
\end{quote}

\noindent where $a_{i}^{(p)},b_{i}^{(p)},g_{s}^{(p)},e_{j}^{(p)}$ are the
polynomials obtained respectively from the polynomials $a_{i},b_{i}%
,g_{s},e_{j}$ in Theorem A, by eliminating those special Schubert classes
$y_{t_{j}}$ with $t_{j}\notin D_{1}(G,p)$ using $q_{j}\cdot a_{j}%
\in\left\langle \omega_{1},\ldots,\omega_{n}\right\rangle $. We obtain

\bigskip

\noindent\textbf{Lemma 3.1.} \textsl{The epimorphism }$\overline{I}_{p}$
\textsl{induces an isomorphism of algebras}

\begin{enumerate}
\item[(3.3)] $H^{\ast}(G/T;\mathbb{F}_{p})=\frac{\mathbb{F}_{p}[\omega
_{1},\cdots,\omega_{n},\text{ }y_{t_{i}}]_{t_{i}\in D_{1}(G,p)}}{\left\langle
a_{i}^{(p)},\text{ }y_{t_{i}}^{r_{i}}-b_{i}^{(p)},\text{ }g_{s}^{(p)},\text{
}e_{j}^{(p)}\right\rangle _{t_{i}\in D_{1}(G,p),t_{s}\in\overline{D}%
_{1}(G,p),1\leq j\leq h}}$\textsl{,}
\end{enumerate}

\noindent\textsl{where, if we set }$\left\langle \omega_{1},\ldots,\omega
_{n}\right\rangle _{\mathbb{F}_{p}}:=\left\langle \omega_{1},\ldots,\omega
_{n}\right\rangle \otimes\mathbb{F}_{p}$\textsl{, then}

\begin{quote}
\textsl{i) }$a_{i}^{(p)},b_{i}^{(p)},g_{s}^{(p)},e_{j}^{(p)}\in\left\langle
\omega_{1},\ldots,\omega_{n}\right\rangle _{\mathbb{F}_{p}}$\textsl{,}

\textsl{ii)} $\#D_{1}(G,p)+\#\overline{D}_{1}(G,p)+h=n$\textsl{.}
\end{quote}

\noindent\textsl{In particular,} $\{a_{i}^{(p)},g_{s}^{(p)},e_{j}^{(p)}%
\}\in\left\langle \omega_{1},\ldots,\omega_{n}\right\rangle _{\mathbb{F}_{p}%
}\cap\ker\overline{I}_{p}$\textsl{.}\hfill$\square$

\bigskip

\noindent\textbf{Remark 3.2. }Let $D(G,p)$ be the degree set of the $n$
polynomials $\{a_{i}^{(p)},g_{s}^{(p)},e_{j}^{(p)}\}$ in (3.3). Since $\deg
a_{i}^{(p)}=t_{i}\in D_{1}(G,p)$ the set\textbf{ }$D(G,p)$ has the partition

\begin{enumerate}
\item[(3.4)] $D(G,p)=D_{1}(G,p)\sqcup D_{2}(G,p)$,
\end{enumerate}

\noindent where $D_{2}(G,p)=D(G,p)\cap2\cdot q(G)$ by iii) of Theorem A.
Consequently, the set $D(G,p)$ consists of precisely $n$ distinct integers.
For examples, we get from the contents of Tables 1 and 2 that

a) If either $G=SU(n)$, $Sp(n)$ or $p\notin\{2,3,5\}$, then

\begin{quote}
$D_{1}(G,p)=\emptyset$ and $D(G,p)=D_{2}(G,p)$ $=2\cdot q(G)$.
\end{quote}

b) If $G=E_{8}$ and $p\in\{2,3,5\}$, then the decomposition (3.4) is

\begin{quote}
$D(E_{8},2)=\left\{  6,10,18,30\right\}  \sqcup\{4,16,24,28\}$;

$D(E_{8},3)=\left\{  8,20\right\}  \sqcup\{4,16,28,36,40,48\}$;

$D(E_{8},5)=\left\{  12\right\}  \sqcup\{4,16,24,28,36,40,48\}$.\hfill
$\square$
\end{quote}

\subsection{The cohomology $H^{\ast}(G;\mathbb{F}_{p})$}

Using the $\mathbb{F}_{p}$-analogues of the maps $[\mathcal{D]}$ and $\kappa$
introduced in Section \S 1.2:

\begin{quote}
$[\mathcal{D}^{\prime}]:\left\langle \omega_{1},\ldots,\omega_{n}\right\rangle
_{\mathbb{F}_{p}}\cap\ker\overline{I}_{p}\rightarrow E_{3}^{\ast
,1}(G;\mathbb{F}_{p}),$

$\kappa^{\prime}:$ $E_{3}^{2k,1}(G;\mathbb{F}_{p})\twoheadrightarrow
E_{\infty}^{2k,1}(G;\mathbb{F}_{p})\subset H^{2k+1}(G;\mathbb{F}_{p})$,
\end{quote}

\noindent we construct, in term of the presentation (3.3) of $H^{\ast
}(G/T;\mathbb{F}_{p})$, a set of generators of the algebra $H^{\ast
}(G;\mathbb{F}_{p})$.

\bigskip

\noindent\textbf{Definition 3.3.} Let $\pi_{p}^{\ast}:H^{\ast}(G/T;\mathbb{F}%
_{p})\rightarrow H^{\ast}(G;\mathbb{F}_{p})$ be the map induced by the bundle
map $\pi$. For a $p$-special Schubert classes $y_{t}$ define the
$\operatorname{mod}p$\textsl{ Schubert cocycles} on $G$ by

\begin{enumerate}
\item[(3.5)] $\overline{x}_{t}:=\pi_{p}^{\ast}(y_{t})\in H^{t}(G;\mathbb{F}%
_{p})$, $t\in D_{1}(G,p)$.
\end{enumerate}

By Remark 3.2, for each $m\in D(G,p)$ there exists a unique polynomial
$c_{m}\in\{a_{i}^{(p)},g_{s}^{(p)},e_{j}^{(p)}\}$ with $\deg c_{m}=m$. Since
$c_{m}\in\left\langle \omega_{1},\ldots,\omega_{n}\right\rangle _{\mathbb{F}%
_{p}}\cap\ker\overline{I}_{p}$ the composition $\kappa^{\prime}\circ
\lbrack\mathcal{D}^{\prime}]$ is applicable to $c_{m}$ to yield the
$\operatorname{mod}p$ cohomology class

\begin{enumerate}
\item[(3.6)] $\zeta_{m-1}:=\kappa^{\prime}\circ\lbrack\mathcal{D}^{\prime
}](c_{m})\in$ $H^{m-1}(G;\mathbb{F}_{p})$,
\end{enumerate}

\noindent to be called a $p$\textsl{-primary class} of $H^{\ast}%
(G;\mathbb{F}_{p})$. Note that if $m\in D_{1}(G,p)$, then $\zeta_{m-1}%
=\theta_{m-1}$ by (1.15).\hfill$\square$

\bigskip

In terms of the $\operatorname{mod}p$ Schubert cocycles $\overline{x}_{t_{i}}%
$, and the $p$-primary classes $\zeta_{m-1}$ just defined, the algebra
$H^{\ast}(G;\mathbb{F}_{p})$ is presented uniformly in the following result.

\bigskip

\noindent\textbf{Theorem 3.4.} \textsl{The inclusion }$\zeta_{m-1}$\textsl{,
}$x_{t}\in H^{\ast}(G;\mathbb{F}_{p})$\textsl{\ induces an isomorphism of
algebras}

\begin{enumerate}
\item[(3.7)] $H^{\ast}(G;\mathbb{F}_{p})=\operatorname{Im}\pi_{p}^{\ast
}\otimes\Delta(\zeta_{m-1})_{m\in D(G,p)}$\textsl{,}
\end{enumerate}

\noindent\textsl{where}

\begin{quote}
\textsl{i) }$\operatorname{Im}\pi_{p}^{\ast}=\mathbb{F}_{p}[\overline{x}%
_{t}]/\left\langle \overline{x}_{t}^{r_{t}}\right\rangle _{t\in D_{1}(G,p)}%
$\textsl{,} $r_{t}=cl(y_{t})$\textsl{,}

\textsl{ii) If }$p\neq2$\textsl{, then} $\Delta(\zeta_{m-1})_{m\in
D(G,p)}=\Lambda(\zeta_{m-1})_{m\in D(G,p)}$\textsl{;}

\textsl{iii) If }$p=2$ \textsl{and} $G$\textsl{ is exceptional, then }%
$\zeta_{m-1}^{2}=0$\textsl{ with the following exceptions}

$\zeta_{3}^{2}=\overline{x}_{6}$\textsl{ for }$G=G_{2},F_{4},E_{6},E_{7}%
,E_{8}$\textsl{;}

$\zeta_{5}^{2}=\overline{x}_{10},\quad\zeta_{9}^{2}=\overline{x}_{18}$\textsl{
for }$G=E_{7},E_{8}$\textsl{;}

$\zeta_{15}^{2}=\overline{x}_{30},$\textsl{ }$\zeta_{23}^{2}=\overline{x}%
_{6}^{6}\overline{x}_{10}$\textsl{ for }$G=E_{8}$\textsl{.}
\end{quote}

\noindent\textbf{Proof.} The proof will make the use of several results
available in the previous paper \cite{DZ1}. In term of the $p$%
\textsl{-transgressive generators} $\{\alpha_{m-1},m\in D(G,p)\}$ of $H^{\ast
}(G;\mathbb{F}_{p})$ introduced in \cite[Definition 3.3]{DZ1}, it has been
shown in \cite[Lemma 3.2]{DZ1} that

\begin{enumerate}
\item[(3.8)] $H^{\ast}(G;\mathbb{F}_{p})=\operatorname{Im}\pi_{p}^{\ast
}\otimes\Delta(\alpha_{m-1})_{m\in D(G,p)}$, where

$\operatorname{Im}\pi_{p}^{\ast}=\mathbb{F}_{p}[\overline{x}_{t}]/\left\langle
\overline{x}_{t}^{r_{t}}\right\rangle _{t\in D_{1}(G,p)}$.
\end{enumerate}

\noindent It has also been shown in \cite[formula (4.2)]{DZ1} that, each
$p$-primary classes $\zeta_{m-1}$ can be expressed as a $\operatorname{Im}%
\pi_{p}^{\ast}$-linear combination of the $p$-transgressive ones with leading
term $\alpha_{m-1}$

\begin{quote}
$\zeta_{m-1}=\alpha_{m-1}+\underset{i\in D(G,p)}{\Sigma}c_{i}\cdot\alpha
_{i-1}$, $c_{i}\in(\operatorname{Im}\pi_{p}^{\ast})^{+}$.
\end{quote}

\noindent Therefore, we obtain (3.7), together with the formula i) of
$\operatorname{Im}\pi_{p}^{\ast}$, from (3.8).

Suppose next that $p\neq2$. With $\deg\zeta_{m-1}\equiv1\operatorname{mod}2$
we get ii) from $\zeta_{m-1}^{2}\equiv0$ by $2\zeta_{m-1}^{2}\equiv
0\operatorname{mod}p$.

To show iii) assume that $G$ is exceptional and $p=2$. By \cite[Theorem
4.1]{DZ1} we have $\zeta_{2s-1}=\alpha_{2s-1}$\ with the following exceptions:

\begin{quote}
$\zeta_{15}=\alpha_{15}+\overline{x}_{6}\alpha_{9}$\textsl{; }$\zeta
_{27}=\alpha_{27}+\overline{x}_{10}\alpha_{17}$\textsl{\ }in\textsl{ }%
$E_{7},E_{8}$\textsl{,}

$\zeta_{23}=\alpha_{23}+\overline{x}_{6}\alpha_{17}$\textsl{\ }in\textsl{
}$E_{7}$\textsl{;}

$\zeta_{23}=\alpha_{23}+\overline{x}_{6}\alpha_{17}+x_{6}^{3}\alpha_{5}%
$\textsl{; }$\quad\zeta_{29}=\alpha_{29}+\overline{x}_{6}^{2}\alpha_{17}%
$\textsl{\ }in\textsl{ }$E_{8}$\textsl{.}
\end{quote}

\noindent The squares $\zeta_{m-1}^{2}$ stated in iii) are verified by the
following relations shown in \cite[(4.8)]{DZ1}

\begin{quote}
$\alpha_{2s-1}^{2}=\left\{
\begin{tabular}
[c]{l}%
$\overline{x}_{6}$ for $s=2$,\\
$\overline{x}_{4s-2}$ for $s=3,5$ and in $E_{7},E_{8}$,\\
$\overline{x}_{30}+\overline{x}_{6}^{2}\overline{x}_{18}$ for $s=8$ and in
$E_{8},$\\
$0$ in the remaining cases.
\end{tabular}
\right.  $.\hfill$\square$
\end{quote}

\subsection{The Bockstein operator $\delta_{p}$ on $H^{\ast}(G;\mathbb{F}%
_{p})$}

In view of the exact sequence (2.5), the\textsl{ }$\operatorname{mod}p$
\textsl{Bockstein operator }on the algebra\textsl{ }$H^{\ast}(G;\mathbb{F}%
_{p})$ is the composition

\begin{quote}
$\delta_{p}=r_{p}\circ\beta_{p}:$ $H^{\ast}(G;\mathbb{F}_{p})\rightarrow
H^{\ast}(G;\mathbb{F}_{p})$
\end{quote}

\noindent which satisfies, as an antiderivation of degree $1$, that

\begin{enumerate}
\item[(3.9)] $\delta_{p}(a\cup r_{p}(c))=\delta_{p}(a)\cup r_{p}(c),a\in
H^{\ast}(G;\mathbb{F}_{p})$, $c\in H^{\ast}(G)$.
\end{enumerate}

\noindent The $\operatorname{mod}p$ Bockstein cohomology $H_{\beta}^{\ast
}(G;\mathbb{F}_{p})$ of $G$ is defined as usual to be the graded quotient
group $\ker\delta_{p}/\operatorname{Im}\delta_{p}$. As preparations to compute
the subgroups $H_{\beta}^{\ast}(G;\mathbb{F}_{p})$ and $\operatorname{Im}%
\delta_{p}$, we clarify the actions of the reduction $r_{p}$ and the
Bockstein\textsl{ }$\beta_{p}$, respectively, in the following two lemmas.

\bigskip

\noindent\textbf{Lemma 3.5.} \textsl{With respect to the presentation (3.7)
of} $H^{\ast}(G;\mathbb{F}_{p})$\textsl{,} \textsl{the reduction }$r_{p}$
\textsl{satisfies that}

\begin{quote}
\textsl{i)} $r_{p}(x_{t_{i}})\equiv\overline{x}_{t_{i}}$ \textsl{if} $t_{i}\in
D_{1}(G,p)$\textsl{, }$0$ \textsl{otherwise;}

\textsl{ii)} $r_{p}(\rho_{l(e_{j})-1})\equiv\zeta_{l(e_{j})-1}$ \textsl{for
}$1\leq j\leq h$\textsl{;}

\textsl{iii)} $r_{p}(\rho_{l(g_{i})-1})\equiv-\overline{x}_{t_{i}}^{r_{i}%
-1}\cdot\theta_{t_{i}-1}$ \textsl{for} $t_{i}\in D_{1}(G,p)$\textsl{;}

\textsl{iv)} $r_{p}(\rho_{l(g_{i})-1})\equiv-q_{i}\cdot\zeta_{l(g_{i})-1}$
\textsl{for} $t_{i}\notin D_{1}(G,p)$\textsl{;}

\textsl{v)} $r_{p}(\mathcal{C}_{I})\equiv\delta_{p}(\theta_{I})$\textsl{ for
}$I\subseteq D_{1}(G,p)$\textsl{,}
\end{quote}

\noindent\textsl{where }$r_{i}=cl(y_{t_{i}}),$ $\theta_{t-1}=\zeta_{t-1}%
$\textsl{, and where }$q_{i}$\textsl{ is co-prime to }$p$\textsl{.}

\bigskip

\noindent\textbf{Proof.} For $t_{i}\in D_{1}(G,p)$ we get $r_{p}(x_{t_{i}%
})=\overline{x}_{t_{i}}$ by $r_{p}\circ\pi^{\ast}=\pi_{p}^{\ast}\circ r_{p}$.
On the other hand, if $t_{i}\notin D_{1}(G,p)$, we find by (3.2) that

\begin{quote}
$r_{p}(x_{t_{i}})\equiv\pi_{p}^{\ast}\circ r_{p}(y_{k_{i}})\equiv\pi_{p}%
^{\ast}(q_{i}a_{i})\equiv0$,
\end{quote}

\noindent where the last equality follows from $a_{i}\in\left\langle
\omega_{1},\ldots,\omega_{n}\right\rangle _{\mathbb{F}_{p}}=\ker\pi_{p}^{\ast
}$. This completes the proof of i).

The relations ii), iii) and iv) will be proven in the same way. In view of

\begin{quote}
$r_{p}\circ\lbrack\mathcal{D}]\equiv\lbrack\mathcal{D}^{\prime}]\circ r_{p}:$
$\left\langle \omega_{1},\ldots,\omega_{n}\right\rangle \cap\ker I\rightarrow
E_{3}^{\ast,1}(G;\mathbb{F}_{p})$
\end{quote}

\noindent the following obvious relations in $E_{2}^{\ast,0}(G;\mathbb{F}%
_{p})=H^{\ast}(G/T;\mathbb{F}_{p})$

\begin{quote}
$r_{p}(e_{j})\equiv e_{j}^{(p)}$ for $1\leq j\leq h$,

$r_{p}(y_{t_{i}}^{r_{i}-1}\cdot f_{i}-p_{i}\cdot g_{i})\equiv-y_{i}^{r_{i}%
-1}\cdot a_{i}^{(p)}$ for $t_{i}\in D_{1}(G,p)$,

$r_{p}(y_{t_{i}}^{r_{i}-1}\cdot f_{i}-p_{i}\cdot g_{i})\equiv-q_{i}\cdot
g_{i}^{(p)}$ for $t_{i}\notin D_{1}(G,p)$
\end{quote}

\noindent imply, respectively, the following relations in $E_{3}^{\ast
,1}(G;\mathbb{F}_{p})$

\begin{quote}
a) $r_{p}[\mathcal{D}](e_{j})\equiv\lbrack\mathcal{D}^{\prime}](e_{j}^{(p)})$,

b) $r_{p}[\mathcal{D}](y_{t_{i}}^{r_{i}-1}\cdot f_{i}-p_{i}\cdot g_{i}%
)\equiv-y_{i}^{r_{i}-1}\cdot\lbrack\mathcal{D}^{\prime}](a_{i}^{(p)})$ if
$t_{i}\in D_{1}(G,p)$,

c) $r_{p}[\mathcal{D}](y_{t_{i}}^{r_{i}-1}\cdot f_{i}-p_{i}\cdot g_{i}%
)\equiv-q_{i}[\mathcal{D}^{\prime}](g_{i}^{(p)})$ if $t_{i}\notin D_{1}(G,p)$.
\end{quote}

\noindent Since the map $\kappa$ is natural with respect to the reduction
$r_{p}$ (i.e. $r_{p}\circ\kappa=\kappa^{\prime}\circ r_{p}$), applying
$\kappa^{\prime}$ to both sides of a), b) and c) we obtain ii), iii) and iv)
from the definition (3.6) of the classes $\zeta_{m-1}$. Note that, in the case
iii), the formula (2.8) is required to derive that

\begin{quote}
$\kappa^{\prime}(y_{i}^{r_{i}-1}\cdot\lbrack\mathcal{D}^{\prime}](a_{i}%
^{(p)})$ $=\overline{x}_{t_{i}}^{r_{i}-1}\cdot\theta_{t_{i}-1}$.
\end{quote}

Finally, the relation v) follows directly from the definition (1.16) of
$\mathcal{C}_{I}$.\hfill$\square$

\bigskip

\noindent\textbf{Lemma 3.6. }\textsl{With respect to the presentation (3.7)
of} $H^{\ast}(G;\mathbb{F}_{p})$\textsl{, the Bockstein homomorphism }%
$\beta_{p}$\textsl{ is given by}

\begin{quote}
\textsl{i)} $\beta_{p}(\zeta_{t_{i}-1})=x_{t_{i}}$ \textsl{for} $t_{i}\in
D_{1}(G,p)$\textsl{;}

\textsl{ii)} $\beta_{p}(\zeta_{m-1})=0$ \textsl{for }$m\in D_{2}%
(G,p)$\textsl{;}

\textsl{iii)} $\beta_{p}(\overline{x}_{t_{i}})=0$ \textsl{for }$t_{i}\in
D_{1}(G,p)$\textsl{.}
\end{quote}

\noindent\textbf{Proof. }For $t_{i}\in D_{1}(G,p)$ the relation $f_{i}=p\cdot
y_{t_{i}}-a_{i}$ on $H^{\ast}(G/T)$ (see in Theorem A) implies the following
the diagram chasing in the short exact sequence (2.6) of $d_{2}$-complexes:

\begin{quote}
$%
\begin{array}
[c]{ccccc}
&  & \mathcal{D}(a_{i}) & \overset{r_{p}}{\rightarrow} & \mathcal{D}^{\prime
}(a_{i}^{(p)})\\
&  & d_{2}\downarrow\quad &  & d_{2}\downarrow\quad\\
y_{t_{i}} & \overset{p}{\longrightarrow} & a_{i} &  & 0
\end{array}
$.
\end{quote}

\noindent It concludes, in the exact sequence (2.7), that $y_{t_{i}}%
=\overline{\beta}_{p}([\mathcal{D}^{\prime}](a_{i}^{(p)}))$ by the convention
$y_{t_{i}}=[y_{t_{i}}\otimes1]$ in Lemma 2.1. Applying $\pi^{\ast}$ to both
sides one gets i) from the relation $\beta_{p}\circ\kappa^{\prime}=\pi^{\ast
}\circ\overline{\beta}_{p}$ by (2.9).

For ii) consider a $p$-primary class $\zeta_{m-1}\in H^{\ast}(G;\mathbb{F}%
_{p})$ with $m\in D_{2}(G,p)$. By the properties ii) and iv) of Lemma 3.4, the
integral primary class $\rho_{m-1}$ satisfies $r_{p}(\rho_{m-1})=q\cdot
\zeta_{m-1}$, where $q$ is co-prime to $p$. One obtains ii) from the exactness
of (2.5).

Finally, with $r_{p}(x_{t_{i}})=\overline{x}_{t_{i}}$ for $t_{i}\in
D_{1}(G,p)$ by i) of Lemma 3.5, one obtain iii) by the exactness of
(2.5).\hfill$\square$

\bigskip

For a subsequence $I=\left\{  i_{1},\cdots,i_{d}\right\}  \subseteq
D_{1}(G,p)$ we put $C_{I}=\delta_{p}(\theta_{I})$. Accordingly, for
$p$-monotone sequences $I,H,L\subseteq D_{1}(G,p)$\textsl{ }introduce in the
polynomial algebra $\mathbb{F}_{p}[\overline{x}_{t_{i}},C_{I}]^{+}$ the
following elements

\begin{quote}
$R_{I}:=(\underset{t_{i}\in I}{\Pi}\overline{x}_{t_{i}}^{r_{i}-1})C_{I}$,
$D_{I}:=\underset{i_{s}\in I}{\Sigma}(-1)^{s-1}y_{i_{s}}\cdot C_{I_{s}}$

$S_{H,L}=C_{H}C_{L}+\underset{h_{s}\in H}{\Sigma}(-1)^{s}\overline{x}_{h_{s}%
}\cdot b_{H_{s}\cap K}\cdot C_{\left\langle H_{s},K\right\rangle }$,
\end{quote}

\noindent where $b_{K}=\underset{t\in K}{\Pi}\zeta_{t-1}^{2}\in
\operatorname{Im}\pi_{p}^{\ast}$ (see (3.12) and (3.13) below).

\bigskip

\noindent\textbf{Theorem 3.7. }\textsl{Let }$G$\textsl{ be a Lie group with
rank }$n$\textsl{, and let }$p$\textsl{ be a prime. Then}

\textsl{i)} $\dim H_{\beta}^{\ast}(G;\mathbb{F}_{p})=2^{n}$\textsl{.}

\textsl{ii) The algebra }$\operatorname{Im}\delta_{p}$ \textsl{has the
presentation}

\begin{enumerate}
\item[(3.10)] $\operatorname{Im}\delta_{p}=\frac{\mathbb{F}_{p}[\overline
{x}_{t_{i}},C_{I}]^{+}}{\left\langle \overline{x}_{t_{i}}^{r_{i}},R_{J}%
,D_{K},S_{H,L}\right\rangle }\otimes\Delta(\zeta_{m-1})_{m\in D_{2}(G,p)}%
$\textsl{,}
\end{enumerate}

\noindent\textsl{where }$t_{i}\in D_{1}(G,p)$\textsl{,} $r_{i}=cl(y_{t_{i}}),$
\textsl{and where} $I,J,K,H,L,\subseteq D_{1}(G,p)$\textsl{ in which
}$\left\vert I\right\vert ,\left\vert J\right\vert ,\left\vert H\right\vert
,\left\vert L\right\vert \geq2$\textsl{, }$\left\vert K\right\vert \geq
3$\textsl{.}

\bigskip

\noindent\textbf{Proof. }According to (3.7) the algebra $H^{\ast}%
(G;\mathbb{F}_{p})$ has the decomposition

\begin{enumerate}
\item[(3.11)] $H^{\ast}(G;\mathbb{F}_{p})=(\operatorname{Im}\pi_{p}^{\ast
}\otimes\Delta(\zeta_{t_{i}-1})_{t_{i}\in D_{1}(G,p)})\otimes\Delta
(\zeta_{m-1})_{m\in D_{2}(G,p)}$
\end{enumerate}

\noindent on which the $\delta_{p}$-action satisfies, by Lemma 3.6, that

\begin{quote}
$\delta_{p}(x)=0$ for $x\in\operatorname{Im}\pi_{p}^{\ast}$,

$\delta_{p}(\zeta_{t_{i}-1})=\overline{x}_{t_{i}}\in\operatorname{Im}\pi
_{p}^{\ast}$ for $t_{i}\in D_{1}(G,p)$,

$\delta_{p}(\zeta_{m-1})=0$ for $m\in D_{2}(G,p)$.
\end{quote}

\noindent These imply that

\begin{quote}
a) The factor $\operatorname{Im}\pi_{p}^{\ast}\otimes\Delta(\zeta_{t_{i}%
-1})_{t_{i}\in D_{1}(G,p)}$ is invariant under $\delta_{p}$, which can be
identified with the Koszul complex $K(\operatorname{Im}\pi_{p}^{\ast})$
associated to the truncated polynomial algebra $\operatorname{Im}\pi_{p}%
^{\ast}$ presented in i) of Theorem 3.4 (see in Section \S 2.2).

b) $\delta_{p}$ acts trivially on the factor $\Delta(\zeta_{m-1})_{m\in
D_{2}(G,p)}$.
\end{quote}

\noindent In particular, let $\delta$ be the restriction of $\delta_{p}$ on
$\operatorname{Im}\pi_{p}^{\ast}\otimes\Delta(\zeta_{t_{i}-1})_{t_{i}\in
D_{1}(G,p)}$. Then

\begin{enumerate}
\item[(3.12)] $\operatorname{Im}\delta_{p}=\operatorname{Im}\delta
\otimes\Delta(\zeta_{m-1})_{m\in D_{2}(G,p)}$.
\end{enumerate}

By a) and b) we get from the K\"{u}nneth formula that

\begin{quote}
$H_{\beta}^{\ast}(G;\mathbb{F}_{p})=H^{\ast}(K(\operatorname{Im}\pi_{p}^{\ast
}))\otimes\Delta(\zeta_{m-1})_{m\in D_{2}(G,p)}$.
\end{quote}

\noindent We obtain i) form $\dim H^{\ast}(K(\operatorname{Im}\pi_{p}^{\ast
}))=2^{\left\vert D_{1}(G,p)\right\vert }$ by (2.12) and $\left\vert
D_{1}(G,p)\right\vert +\left\vert D_{2}(G,p)\right\vert =n$ by Remark 3.2.

For ii) we note that the cup product on $H^{\ast}(G;\mathbb{F}_{p})$ furnishes
the complex $K(\operatorname{Im}\pi_{p}^{\ast})$ with the canonical algebra
structure $b=\{\theta_{m-1}^{2}$, $m\in D_{1}(G,p)\}$, where

\begin{enumerate}
\item[(3.13)] $\theta_{m-1}^{2}=0\in\operatorname{Im}\pi_{p}^{\ast}$ if
$p\equiv1\operatorname{mod}2$ (see Remark 2.6);

\item[(3.14)] $\theta_{m-1}^{2}\in\operatorname{Im}\pi_{2}^{\ast}$ if $p=2$
(by the $\operatorname{mod}2$ analogue of Lemma 2.2).
\end{enumerate}

\noindent Thus, we can derive (3.10) from (3.12) and (2.14).\hfill$\square$

\section{The structure of the ring $H^{\ast}(G)$}

In terms of the cohomology classes $x_{t_{i}}$, $\rho_{2l-1}$ and
$\mathcal{C}_{I}$ constructed in Section \S 1.2, we establish in Theorems 4.1,
4.4 and 4.7 general expressions of the three basic components of the
cohomology $H^{\ast}(G)$: the subring $\operatorname{Im}\pi^{\ast}$, the free
part $\mathcal{F}(G)$ and the torsion ideal $\tau_{p}(G)$, respectively. These
results are summarized in Section \S 4.4 to give a proof of Theorem B.

\subsection{The subring $\operatorname{Im}\pi^{\ast}\subset H^{\ast}(G)$}

For a prime $p$ the proof of Lemma 2.1 is valid to show from (3.3) that

\begin{quote}
$E_{3}^{\ast,0}(G;\mathbb{F}_{p})=\mathbb{F}_{p}[y_{t_{i}}]/\left\langle
y_{t_{i}}^{r_{i}}\right\rangle _{i\in D_{1}(G,p)}$.
\end{quote}

\noindent Moreover, with $\pi_{p}^{\ast}(y_{t_{i}})=\overline{x}_{t_{i}}$ by
(3.5) we get by formula (3.7) that the map $\pi_{p}^{\ast}$ induces an isomorphism

\begin{quote}
$E_{3}^{\ast,0}(G;\mathbb{F}_{p})\overset{\cong}{\rightarrow}\operatorname{Im}%
\pi_{p}^{\ast}=\mathbb{F}_{p}[\overline{x}_{t_{i}}]/\left\langle \overline
{x}_{t_{i}}^{r_{t_{i}}}\right\rangle _{i\in D_{1}(G,p)}$.
\end{quote}

\noindent Thus, we can deduce from Lemma 2.1, as well as the relation $\pi
_{p}^{\ast}\circ r_{p}=r_{p}\circ\pi^{\ast}$ on $E_{3}^{\ast,0}(G)$, the
following result.

\bigskip

\noindent\textbf{Theorem 4.1. }\textsl{The map }$\pi^{\ast}:E_{3}^{\ast
,0}(G)\rightarrow\operatorname{Im}\pi^{\ast}$ \textsl{in (2.2) is an
isomorphism.}

\textsl{In particular, if }$\left\{  y_{t_{1}},\cdots,y_{t_{k}}\right\}
$\textsl{ is a set of special Schubert classes on }$G/T$\textsl{,}

\begin{enumerate}
\item[(4.1)] $\operatorname{Im}\pi^{\ast}=\mathbb{Z}[x_{t_{1}},\cdots
,x_{t_{k}}]/\left\langle p_{i}x_{t_{i}},x_{t_{i}}^{r_{i}}\right\rangle _{1\leq
i\leq k}$
\end{enumerate}

\noindent\textsl{where} $p_{i}=tor(y_{t_{i}})$\textsl{,} $r_{i}=cl(y_{t_{i}}%
)$\textsl{.}\hfill$\square$

\bigskip

It is straightforward to see from (4.1) that

\bigskip

\noindent\textbf{Corollary 4.2. }\textsl{The reduction }$r_{p}:H^{\ast
}(G)\rightarrow H^{\ast}(G;\mathbb{F}_{p})$\textsl{ restricts to an
isomorphism}

\begin{enumerate}
\item[(4.2)] $r_{p}:\operatorname{Im}\pi^{\ast}\cap\tau_{p}(G)\rightarrow
\operatorname{Im}(\pi_{p}^{\ast})^{+}$.\hfill$\square$
\end{enumerate}

\bigskip

\noindent\textbf{Example 4.3.} In \cite{G} Grothendieck defined \textsl{Chow
ring} of the reductive algebraic group $G^{c}$ corresponding to $G$ to be the subring

\begin{quote}
$\mathcal{A}(G^{c}):=\operatorname{Im}\left\{  \pi^{\ast}:H^{\ast
}(G/T)\rightarrow H^{\ast}(G)\right\}  $.
\end{quote}

\noindent Since the Chow ring $\mathcal{A}(G/T)$ of the projective variety
$G/T$ is canonically isomorphic to the integral cohomology $H^{\ast}(G/T)$,
$\mathcal{A}(G^{c})$ consists of those elements of $H^{\ast}(G)$ that can be
realized by the algebraic cocycles on $G/T$ via $\pi$.

In view of formula (4.1), the degree set $\left\{  t_{1},\cdots,t_{k}\right\}
$ of the special Schubert classes on $G/T$, as well as the torsion indices
$tor(y_{t_{i}})$ and the cup lengths $cl(y_{t_{i}})$, are the numerical
invariants of $G$ required to express $\mathcal{A}(G^{c})$ algebraically. In
particular, inputting the values of these invariants given by Table 2, we get
from (4.1) the following expressions of the ring $\mathcal{A}(G^{c})$ in terms
of the Schubert cocycles on $G$.

\begin{center}
{\normalsize Table 3. The Chow rings of exceptional Lie groups}%

\begin{tabular}
[c]{l||l}\hline
$G$ & $\mathcal{A}(G^{c})=\operatorname{Im}\pi^{\ast}\{H^{\ast}%
(G/T)\rightarrow H^{\ast}(G)\}$\\\hline\hline
$G_{2}$ & $\frac{\mathbb{Z}[x_{6}]}{\left\langle 2x_{6},x_{6}^{2}\right\rangle
}$\\\hline
$F_{4}$ & $\frac{\mathbb{Z}[x_{6},x_{8}]}{\left\langle 2x_{6},x_{6}^{2}%
,3x_{8},x_{8}^{3}\right\rangle }$\\\hline
$E_{6}$ & $\frac{\mathbb{Z}[x_{6},x_{8}]}{\left\langle 2x_{6},x_{6}^{2}%
,3x_{8},x_{8}^{3}\right\rangle }$\\\hline
$E_{7}$ & $\frac{\mathbb{Z}[x_{6},x_{8},x_{10},x_{18}]}{\left\langle
2x_{6},3x_{8},2x_{10},2x_{18},x_{6}^{2},x_{8}^{3},x_{10}^{2},x_{18}%
^{2}\right\rangle }$\\\hline
$E_{8}$ & $\frac{\mathbb{Z}[x_{6},x_{8},x_{10},x_{12},x_{18},x_{20},x_{30}%
]}{\left\langle 2x_{6},3x_{8},2x_{10},5x_{12},2x_{18},3x_{20},2x_{30}%
,x_{6}^{8},x_{8}^{3},x_{10}^{4},x_{12}^{5},x_{18}^{2},x_{20}^{3},x_{30}%
^{2}\right\rangle }$\\\hline
\end{tabular}
.
\end{center}

We note that presentations of the Chow rings $\mathcal{A}(G^{c})$ for
$G=Spin(n)$, $SO(n),G_{2}$ and $F_{4}$, in term of Schubert classes on $G/T$,
have been obtained by Marlin in \cite{M2}.\hfill$\square$

\subsection{The free part $\mathcal{F}(G)\subseteq H^{\ast}(G)$}

According to Example 1.14, if $G$ is a Lie group with rank $n$, the totality
of the primary classes of $G$ is $\{\rho_{2l-1}\in H^{\ast}(G),l\in q(G)\}$,
where $q(G)=\{l_{1},\cdots,l_{n}\}$ is the degree set of basic homogeneous
$W$-invariants of $G$ \cite{K}.

\bigskip

\noindent\textbf{Theorem 4.4. }$\mathcal{F}(G)=\Delta(\rho_{2l_{1}-1}%
,\cdots,\rho_{2l_{n}-1}).$

\bigskip

\noindent\textbf{Proof. }Recall that the group $\Delta(\rho_{2l_{1}-1}%
,\cdots,\rho_{2l_{n}-1})$ has the basis consisting of all the square free
monomials in $\rho_{2l_{1}-1},\cdots,\rho_{2l_{n}-1}$:

\begin{quote}
$\Phi=\{1,\rho_{I}=\rho_{2l_{i_{1}}-1}\cdots\rho_{2l_{i_{k}}-1}\}$,
$I=\{i_{1},\ldots,i_{k}\}\subseteq\{1,2,\ldots,n\}$.
\end{quote}

\noindent It suffices to show that $\Phi$ is a basis of $\mathcal{F}(G)$.
Suppose that $m=\dim G$.

By (1.17) the full product $\rho=\rho_{2l_{1}-1}\cdots\rho_{2l_{n}-1}\in\Phi$
belongs to the top degree cohomology group $H^{m}(G)=\mathbb{Z}$. On the other
hand, by formula (3.7), for any prime $p$ the top degree cohomology group
$H^{m}(G;\mathbb{F}_{p})=\mathbb{F}_{p}$ is generated by the product

\begin{quote}
$\omega_{p}=\underset{t_{i}\in D_{1}(G,p)}{\Pi}\overline{x}_{t_{i}}^{r_{t_{i}%
}-1}\underset{m\in D(G,p)}{\Pi}\zeta_{m-1}$.
\end{quote}

\noindent From properties ii), iii) and iv) of Lemma 3.5 one finds that

\begin{quote}
$r_{p}(\rho)=r_{p}(\rho_{2l_{1}-1})\cdot\cdots\cdot r_{p}(\rho_{2l_{n}%
-1})\equiv q\cdot\omega_{p}$,
\end{quote}

\noindent where $q$ is an integer co-prime to $p$. Therefore, $\rho$ must be a
generator of the group $H^{m}(G)=\mathbb{Z}$. In particular, the set $\Phi$ of
$2^{n}$ monomials is linearly independent in $H^{\ast}(G)$. Since $\dim
H^{\ast}(G)\otimes\mathbb{Q=}2^{n}$ it suffices to show that $\Phi$ spans a
direct summand of $H^{\ast}(G)$.

Assume, on the contrary, that there exist a monomial $\rho_{I}\in\Phi$, a
class $\xi\in H^{\ast}(G)$, as well as some integer $a>1$, so that a relation
of the form $\rho_{I}=a\cdot\xi$ holds in $H^{\ast}(G)$. Letting $\overline
{I}$ be the complement of $I\subseteq\{1,\cdots,n\}$, and multiplying both
sides by the class $\rho_{\overline{I}}$ yield

\begin{quote}
$\rho=(-1)^{r}a\cdot(\xi\cup\rho_{\overline{I}})$ for some $r\in\mathbb{Z}$.
\end{quote}

\noindent This contradicts to the fact that $\rho$ generates the group
$H^{m}(G)=\mathbb{Z}$. This completes the proof.\hfill$\square$

\bigskip

\noindent\textbf{Example 4.5.} As a supplement to Theorem 4.3 we show that, if
$G$ is exceptional, then

\begin{quote}
$\mathcal{F}(G_{2})=\Delta(\varrho_{3})\otimes\Lambda(\varrho_{11})$

$\mathcal{F}(F_{4})=\Delta(\varrho_{3})\otimes\Lambda(\varrho_{11}%
,\varrho_{15},\varrho_{23})$

$\mathcal{F}(E_{6})=\Delta(\varrho_{3})\otimes\Lambda(\varrho_{9},\varrho
_{11},\varrho_{15},\varrho_{17},\varrho_{23})$

$\mathcal{F}(E_{7})=\Delta(\varrho_{3})\otimes\Lambda(\varrho_{11}%
,\varrho_{15},\varrho_{19},\varrho_{23},\varrho_{27},\varrho_{35})$

$\mathcal{F}(E_{8})=\Delta(\varrho_{3},\varrho_{15},\varrho_{23}%
)\otimes\Lambda(\varrho_{27},\varrho_{35},\varrho_{39},\varrho_{47}%
,\varrho_{59})$
\end{quote}

\noindent where

\begin{quote}
$\rho_{3}^{2}=x_{6}$ for $G=G_{2},F_{4},E_{6},E_{7},E_{8}$;

$\rho_{15}^{2}=x_{30}$, $\rho_{23}^{2}=x_{6}^{6}x_{10}$ for $G=E_{8}$.
\end{quote}

We show the results for the relatively nontrivial case $G=E_{8}$, for which
one has, by the contents in the last column of Table 2, that

\begin{quote}
$\tau_{2}(E_{8})\cap\operatorname{Im}\pi^{\ast}=\frac{\mathbb{F}_{2}%
[x_{6},x_{10},x_{18},x_{30}]^{+}}{\left\langle x_{6}^{8},x_{10}^{4},x_{18}%
^{2},x_{30}^{2}\right\rangle }$.
\end{quote}

\noindent Since $\rho_{2l-1}^{2}\in\tau_{2}(E_{8})\cap\operatorname{Im}%
\pi^{\ast}$ by Lemma 2.2 and since, by Corollary 4.2, $r_{2}$ maps $\tau
_{2}(E_{8})\cap\operatorname{Im}\pi^{\ast}$ isomorphically onto

\begin{quote}
$(\operatorname{Im}\pi_{2}^{\ast})^{+}=$ $\frac{\mathbb{F}_{2}[\overline
{x}_{6},\overline{x}_{10},\overline{x}_{18},\overline{x}_{30}]^{+}%
}{\left\langle \overline{x}_{6}^{8},\overline{x}_{10}^{4},\overline{x}%
_{18}^{2},\overline{x}_{30}^{2}\right\rangle }$,
\end{quote}

\noindent the squares $\rho_{2l-1}^{2}$ are determined by $r_{2}(\rho
_{2l-1})^{2}$. The proof is completed by the results tabulated below

\begin{center}%
\begin{tabular}
[c]{l||l|l|l|l|l|l|l|l}%
$\rho_{2l-1}$ & $\rho_{3}$ & $\rho_{15}$ & $\rho_{23}$ & $\rho_{27}$ &
$\rho_{35}$ & $\rho_{39}$ & $\rho_{47}$ & $\rho_{59}$\\\hline\hline
$r_{2}(\rho_{2l-1})$ & $\zeta_{3}$ & $\zeta_{15}$ & $\zeta_{23}$ & $\zeta
_{27}$ & $\overline{x}_{18}\zeta_{17}$ & $\overline{x}_{10}^{3}\zeta_{9}$ &
$\overline{x}_{6}^{7}\zeta_{5}$ & $\overline{x}_{30}\zeta_{29}$\\\hline
$r_{2}(\rho_{2l-1})^{2}$ & $\overline{x}_{6}$ & $\overline{x}_{15}$ &
$\overline{x}_{6}^{6}\overline{x}_{10}$ & $0$ & $0$ & $0$ & $0$ & $0$\\\hline
\end{tabular}
,
\end{center}

\noindent where the contents in the second row follow from properties ii),
iii) and iv) of Lemma 3.5, and where the results in the third row have been
shown in iii) of Theorem 3.4.\hfill$\square$\hfill

\subsection{The torsion ideal $\tau_{p}(G)\subset H^{\ast}(G)$}

For a vector space $V$ graded by finite dimensional subspaces $V^{k}$
($k\geq0$) denote its Poincar\'{e} polynomial by $P_{t}(V)=\Sigma\dim
V^{k}\cdot t^{k}$. For a finite CW-complex $X$ and a prime $p$, let
$\delta_{p}=r_{p}\circ\beta_{p}$ be the Bockstein operator on $H^{\ast
}(X;\mathbb{F}_{p})$. From the short exact sequence

\begin{quote}
$0\rightarrow\ker\delta_{p}\rightarrow H^{\ast}(X;\mathbb{F}_{p}%
)\overset{\delta_{p}}{\rightarrow}\operatorname{Im}\delta_{p}\rightarrow0$
\end{quote}

\noindent together with $\ker\delta_{p}=\operatorname{Im}\delta_{p}\oplus
H_{\beta}^{\ast}(X;\mathbb{F}_{p})$, we get the following (e.g. \cite[Section
30.4]{BH1})

\begin{quote}
$P_{t}(H^{\ast}(X;\mathbb{F}_{p}))=P_{t}(H_{\beta}^{\ast}(X;\mathbb{F}%
_{p}))+(1+t^{-1})P_{t}(\operatorname{Im}\delta_{p})$.
\end{quote}

\noindent We apply this to show that

\bigskip

\noindent\textbf{Lemma 4.6. }If $\dim H^{\ast}(X;\mathbb{Q})=\dim H_{\beta
}^{\ast}(X;\mathbb{F}_{p})$, then the reduction $r_{p}:H^{\ast}(X)\rightarrow
H^{\ast}(X;\mathbb{F}_{p})$ restricts to an isomorphism

\begin{enumerate}
\item[(4.3)] $r_{p}:\tau_{p}(X)\overset{\cong}{\rightarrow}\operatorname{Im}%
\delta_{p}$.
\end{enumerate}

\noindent\textbf{Proof.} With respect to the decomposition $H^{\ast
}(X)=\mathcal{F}(X)\oplus_{p}\tau_{p}(X)$ by (1.18) the reduction $r_{p}$
induces a monomorphism

\begin{quote}
$r_{p}\otimes1:$ $\mathcal{F}(X)\otimes\mathbb{F}_{p}=\mathcal{F}%
(X)/p\cdot\mathcal{F}(X)\rightarrow H_{\beta}^{\ast}(X;\mathbb{F}%
_{p})\subseteq\ker\delta_{p}$.
\end{quote}

\noindent It implies, if

\begin{quote}
$\dim H^{\ast}(X;\mathbb{Q})(=\mathcal{F}(X)\otimes\mathbb{Q})=\dim H_{\beta
}^{\ast}(X;\mathbb{F}_{p})$,
\end{quote}

\noindent then $r_{p}\otimes1$ is an isomorphism. In particular, we get from

\begin{quote}
$P_{t}(H^{\ast}(X;\mathbb{Q}))=P_{t}(\mathcal{F}(X)\otimes\mathbb{F}%
_{p})=P_{t}(H_{\beta}^{\ast}(X;\mathbb{F}_{p}))$
\end{quote}

\noindent the relation

\begin{quote}
$P_{t}(H^{\ast}(X;\mathbb{F}_{p}))=P_{t}(H^{\ast}(X;\mathbb{Q}))+(1+t^{-1}%
)P_{t}(\operatorname{Im}\delta_{p})$.
\end{quote}

\noindent Lemma 4.6 follows now from Borel-Hirzebruch \cite[Section 30.4]%
{BH1}.\hfill$\square$

\bigskip

Suppose now that $G$ is a Lie group and $p$ is a prime. By Theorem 4.1, for
each $t_{i}\in D_{1}(G,p)$ the Schubert cocycle $x_{t_{i}}=\pi^{\ast}%
(y_{t_{i}})$ satisfies $x_{t_{i}}\in\tau_{p}(G)$. As in (1.16), for each
$p$-monotone sequence $I\subseteq D_{1}(G,p)$ we put

\begin{quote}
$\mathcal{C}_{I}=\beta_{p}(\theta_{I})\in\tau_{p}(G)$.
\end{quote}

\noindent Accordingly, for subsequences $I,H,L\subseteq D_{1}(G,p)$ consider
in the polynomial algebra $\mathbb{F}_{p}[x_{t_{i}},\mathcal{C}_{I}]_{t_{i}\in
D_{1}(G,p),\text{ }I\subseteq D_{1}(G,p)}^{+}$ the following elements

\begin{enumerate}
\item[(4.4)] $\mathcal{R}_{I}:=(\underset{t_{i}\in I}{\Pi}x_{t_{i}}^{r_{i}%
-1})\mathcal{C}_{I}$, $\mathcal{D}_{I}:=\underset{i_{s}\in I}{\Sigma
}(-1)^{s-1}y_{i_{s}}\cdot C_{I_{s}}$,

\item[(4.5)] $\mathcal{S}_{H,L}:=\mathcal{C}_{H}\mathcal{C}_{L}+\underset
{h_{s}\in H=\left\{  h_{1},\cdots,h_{d}\right\}  }{\Sigma}(-1)^{s}x_{h_{s}%
}\cdot b_{H_{s}\cap L}^{\ast}\cdot\mathcal{C}_{\left\langle H_{s}%
,L\right\rangle }$,
\end{enumerate}

\noindent where $r_{j}=cl(y_{t_{j}})$, and where $b_{K}^{\ast}\in
\operatorname{Im}\pi^{\ast}\cap\tau_{p}(G)$ is the unique element that
satisfies by Corollary 4.2 the following relation (see (3.13) and (3.14))

\begin{quote}
$r_{p}(b_{K}^{\ast})=\underset{t\in K}{\Pi}\zeta_{t-1}^{2}\in
(\operatorname{Im}\pi_{p}^{\ast})^{+}$.
\end{quote}

\bigskip

\noindent\textbf{Theorem 4.7. }\textsl{For a Lie group }$G$\textsl{ and prime
}$p$\textsl{ we have}

\begin{enumerate}
\item[(4.6)] $\tau_{p}(G)=\frac{\mathbb{F}_{p}[x_{t_{i}},\mathcal{C}_{I}]^{+}%
}{\left\langle x_{t_{i}}^{r_{i}},\mathcal{R}_{J},\mathcal{D}_{K}%
,\mathcal{S}_{H,L}\right\rangle }\otimes\Delta(\rho_{m-1})_{m\in D_{2}(G,p)}%
$\textsl{,}
\end{enumerate}

\noindent\textsl{where }$t_{i}\in D_{1}(G,p)$\textsl{, and where
}$I,J,K,H,L\subseteq D_{1}(G,p)$\textsl{ in which }$\left\vert I\right\vert
,\left\vert J\right\vert ,\left\vert H\right\vert ,\left\vert L\right\vert
\geq2$\textsl{, }$\left\vert K\right\vert \geq3$\textsl{.}

\bigskip

\noindent\textbf{Proof.} Assume that the rank of the group $G$ is $n$. Since

\begin{quote}
$\dim H^{\ast}(G;\mathbb{Q})=\dim H_{\beta}^{\ast}(G;\mathbb{F}_{p})=2^{n}$,
\end{quote}

\noindent by Theorems 3.7 and 4.4, the restriction $r_{p}:\tau_{p}%
(G)\rightarrow\operatorname{Im}\delta_{p}$ of $r_{p}$ on $\tau_{p}(G)$ is an
isomorphism by Lemma 4.6. With

\begin{quote}
i) $r_{p}(\mathcal{C}_{I})\equiv\delta_{p}(\theta_{I})=C_{I}$,

ii) $r_{p}(x_{t})\equiv\overline{x}_{t}$ for $t\in D_{1}(G,p)$ (by i) of Lemma
3.5), and

iii) $r_{p}(\rho_{m-1})\equiv q_{m}\cdot\zeta_{m-1}$ for $m\in D_{2}(G,p)$ (by
ii), iv) of Lemma 3.5), where $q_{m}$ is co-prime to $p$,
\end{quote}

\noindent one translates the formula (3.10) of $\operatorname{Im}\delta_{p}$
to the desired formula (4.6) of $\tau_{p}(X)$.\hfill$\square$

\bigskip

\noindent\textbf{Example 4.8.} Let $G$ be an exceptional Lie group. Applying
Theorem 4.6 we justify the expressions of the torsion ideals $\tau_{p}(G)$
stated in Theorem B.

Since $D_{1}(G,p)=\emptyset$ implies $\tau_{p}(G)=0$ by Theorem 4.6, we can
assume below that $D_{1}(G,p)\neq\emptyset$. For convenience the calculation
is divided into three cases.

\textbf{Case 1.} If $(G,p)=(G_{2},2),(F_{4},2),(F_{4},3),(E_{6},2),(E_{6}%
,3),(E_{7},3)$ or $(E_{8},5)$ then the set $D_{1}(G,p)$ is a singleton $\{t\}$
by Table 2. It follows that the generators $\mathcal{C}_{I}\in\tau_{p}(G)$
with $\left\vert I\right\vert \geq2$ are absent, hence the formula (4.6) turns
to be

\begin{quote}
$\tau_{p}(G)=\frac{\mathbb{F}_{p}[x_{t}]^{+}}{\left\langle x_{t}%
^{r}\right\rangle }\otimes\Delta(\rho_{m-1})_{m\in D_{2}(G,p)}$, $r=cl(y_{t})$.
\end{quote}

\noindent In addition, if either $p=3$ or $5$, we can replace the factor
$\Delta(\rho_{m-1})_{m\in D_{2}(G,p)}$ by the exterior algebra $\Lambda
(\rho_{m-1})_{m\in D_{2}(G,p)}$, because $\deg(\rho_{m-1})\equiv
1\operatorname{mod}2$ implies that $\rho_{m-1}^{2}\equiv0\operatorname{mod}p$.
Finally, if $p=2$, the squares $\rho_{m-1}^{2}$ with $m\in D_{2}(G,p)$ have
been evaluated in Example 4.5. This justified, for the current situations, the
presentations of the ideal $\tau_{p}(G)$ stated in Theorem B.

\textbf{Case 2.} If $(G,p)=(E_{8},3)$ then $D_{1}(G,p)=\{8,20\}$ by Table 2.
The formula (4.6) for $\tau_{3}(E_{8})$ turns to be

\begin{center}
$\frac{\mathbb{F}_{3}[x_{8},x_{20},\mathcal{C}_{\{8,20\}}]^{+}}{\left\langle
x_{8}^{3},x_{20}^{3},\mathcal{R}_{\{8,20\}},\mathcal{S}_{\{8,20\},\{8,20\}}%
\right\rangle }\otimes\Delta(\varrho_{3},\varrho_{15},\varrho_{27}%
,\varrho_{35},\varrho_{39},\varrho_{47})$,
\end{center}

\noindent where $\mathcal{R}_{\{8,20\}}=x_{8}^{2}x_{20}^{2}\mathcal{C}%
_{\{8,20\}}$ by the definition of $\mathcal{R}_{J}$, $\mathcal{S}%
_{\{8,20\},\{8,20\}}=\mathcal{C}_{\{8,20\}}^{2}$ by (4.5), and where
$\varrho_{2l-1}^{2}\equiv0\operatorname{mod}3$ because of $2\varrho_{2l-1}%
^{2}=0$ by Lemma 2.2. These verifies the formula of $\tau_{3}(E_{8})$ stated
in v) of Theorem B.

\textbf{Case 3. }For the remaining cases $(G,p)=(E_{7},2)$ or $(E_{8},2)$ one
reads from Table 2 that

\begin{quote}
$D_{1}(E_{7},2)=\{6,10,18\}$ and $D_{1}(E_{8},2)=\{6,10,18,30\}$,
\end{quote}

\noindent indicating that the generators $\mathcal{C}_{I}$ of $\tau_{2}(G)$,
as well as the relations $\mathcal{R}_{J},\mathcal{D}_{K}$, $\mathcal{S}%
_{H,L}$ on $\tau_{2}(G)$, are too many to be stated explicitly. Instead, we
point out that, for any subsequences $H,L\subseteq D_{1}(G,2)$, formula (4.5)
is practical to express the relations $\mathcal{S}_{H,L}$ on $\tau_{2}(G)$ as
explicit polynomials in the generators $x_{t_{i}}$ and $\mathcal{C}_{I}$. For
example, the calculations in $H^{\ast}(E_{8},\mathbb{F}_{2})$

\begin{quote}
$C_{\{6,10\}}C_{\{6,10\}}=\overline{x}_{6}\zeta_{9}^{2}C_{\{6\}}+\overline
{x}_{10}\zeta_{5}^{2}C_{\{10\}}=\overline{x}_{6}^{2}\overline{x}%
_{18}+\overline{x}_{10}^{3}$,

$C_{(6,10)}C_{(6,18)}=\overline{x}_{10}\zeta_{5}^{2}C_{\{18\}}+\overline
{x}_{6}C_{\{6,10,18\}}=\overline{x}_{10}^{2}\overline{x}_{18}+\overline{x}%
_{6}C_{\{6,10,18\}}$
\end{quote}

\noindent correspond, under the isomorphism (4.3), to the following relations
on $\tau_{2}(E_{8})$

\begin{quote}
$\mathcal{S}_{\{6,10\},\{6,10\}}=\mathcal{C}_{\{6,10\}}\mathcal{C}%
_{\{6,10\}}+x_{6}^{2}x_{18}+x_{10}^{3}$;

$\mathcal{S}_{\{6,10\},\{6,18\}}=\mathcal{C}_{(6,10)}\mathcal{C}%
_{(6,18)}+x_{10}^{2}x_{18}+x_{6}\mathcal{C}_{\{6,10,18\}}$.
\end{quote}

\noindent These computations show that the relations of the type
$\mathcal{S}_{H,L}$ may be nontrivial.

Finally, we remark that, if $G=E_{7}$, there is only one relation of the type
$\mathcal{D}_{K}$ on \hfill$\tau_{2}(E_{7})$ due to the constraints
$\left\vert K\right\vert \geq3$ and $K\subseteq\{6,10,18\}$, which by (4.4) is

\begin{quote}
$\mathcal{D}_{\{6,10,18\}}=x_{6}\mathcal{C}_{\{10,18\}}-x_{10}\mathcal{C}%
_{\{6,18\}}+x_{18}\mathcal{C}_{\{6,10\}}$.\hfill$\square$
\end{quote}

\subsection{The action of the free part $\mathcal{F}(G)$ on the ideal
$\tau_{p}(G)$}

Since the subgroup $\tau_{p}(G)$ of $H^{\ast}(G)$ is an ideal, the cup product
on $H^{\ast}(G)$ defines an action of the free part $\mathcal{F}(G)$ on
$\tau_{p}(G)$:

\begin{quote}
$\mathcal{F}(G)\times\tau_{p}(G)\rightarrow\tau_{p}(G)$.
\end{quote}

\noindent In views of expressions of $\mathcal{F}(G)$ and $\tau_{p}(G)$ given
respectively by Theorems 4.4 and 4.7, to clarifies this action it suffices to
find a formula $\mathcal{H}_{i,I}$ that expresses each product $\rho
_{l(g_{i})-1}\cdot\mathcal{C}_{K}$, with $t_{i}\in D_{1}(G,p)$ and $K\subseteq
D_{1}(G,p)$, as an element of $\tau_{p}(G)$. Here, if $K=\{t\}$ is a
singleton, then $\mathcal{C}_{K}:=\beta_{p}(\theta_{t-1})=x_{t}$ by i) of
Lemma 3.6.

\bigskip

\noindent\textbf{Theorem 4.9.} \textsl{For each }$t_{i}\in D_{1}(G,p)$\textsl{
and }$K\subseteq D_{1}(G,p)$\textsl{ the relation }$\mathcal{H}_{i,K}$
\textsl{is given by the three possibilities}

\begin{enumerate}
\item[(4.7)] $\rho_{l(g_{i})-1}\cdot\mathcal{C}_{K}=\left\{
\begin{tabular}
[c]{l}%
$x_{t_{i}}^{r_{i}-1}\cdot\mathcal{C}_{K\cup\{t_{i}\}}$ \textsl{if}
$t_{i}\notin K$;\\
$0$ \textsl{if }$p$ \textsl{is odd and }$t_{i}\in K$;\\
$x_{t_{i}}^{r_{i}-1}\cdot(\theta_{t_{i}-1}^{2})^{\ast}\cdot\mathcal{C}%
_{K_{t_{i}}}$ \textsl{if} $p=2$ \textsl{and} $t_{i}\in K$\textsl{,}%
\end{tabular}
\ \ \ \right.  $
\end{enumerate}

\noindent\textsl{where }$K_{t_{i}}$ \textsl{is the complement of }$t_{i}\in
K$\textsl{, and where} $(\theta_{t_{i}-1}^{2})^{\ast}\in\operatorname{Im}%
\pi^{\ast}\cap\tau_{2}(G)$ \textsl{denotes the unique element that satisfies,
under the isomorphism (4.2), that}

\begin{quote}
$r_{2}(\theta_{t_{i}-1}^{2})^{\ast}=\theta_{t_{i}-1}^{2}\in\operatorname{Im}%
(\pi_{2}^{\ast})^{+}$\textsl{.}
\end{quote}

\noindent\textbf{Proof.} Since the reduction $r_{p}$ restricts to an
isomorphism $\tau_{p}(G)\rightarrow\operatorname{Im}\delta_{p}$ by (4.6), we
can find the expression of the product $\rho_{l(g_{i})-1}\cdot\mathcal{C}%
_{K}\in\tau_{p}(G)$ by computing with its $r_{p}$ image in $H^{\ast
}(G;\mathbb{F}_{p})$. Precisely, from $C_{K}=\delta_{p}(\theta_{K})$ and
$r_{p}(\rho_{l(g_{i})-1})\equiv-\overline{x}_{t_{i}}^{r_{\overline{t}}-1}%
\cdot\theta_{t_{i}-1}$ (by iii) of Lemma 3.5) one finds that

\begin{quote}
$r_{p}(\rho_{l(g_{i})-1}\cdot\mathcal{C}_{K})\equiv-\overline{x}_{t_{i}%
}^{r_{\overline{t}}-1}\delta_{p}(\theta_{t_{i}-1}\theta_{K})$.
\end{quote}

\noindent The proof of (4.7) is completed by

\begin{quote}
$\theta_{t_{i}-1}\theta_{K}\equiv\left\{
\begin{tabular}
[c]{l}%
$\theta_{K\cup\{t_{i}\}}\text{ if }t_{i}\notin K\text{;\quad}$\\
$0\text{ if}$ $p$ is odd and $t_{i}\in K\text{;}$\\
$\theta_{t_{i}-1}^{2}\theta_{K_{t_{i}}}\text{ if }p=2\text{ and }t_{i}\in K$.
\end{tabular}
\ \ \right.  $\hfill$\square$
\end{quote}

\bigskip

We are ready to show Theorem B stated in Section \S 1.3.

\bigskip

\noindent\textbf{Proof of Theorem B. }In Theorem B, the expressions of the
cohomology $H^{\ast}(G)$ in the form

\begin{quote}
$H^{\ast}(G)=\mathcal{F}(G)\oplus_{p}\tau_{p}(G)$,
\end{quote}

\noindent together with the presentations of the summands $\mathcal{F}(G)$ and
$\tau_{p}(G)$, have been shown by Examples 4.5 and 4.8. It remains to add
those relations $\mathcal{H}_{i,K}$ that describe the action of $\mathcal{F}%
(G)$ on $\tau_{p}(G)$, where $t_{i}\in D_{1}(G,p)$ and $K\subseteq D_{1}%
(G,p)$. To this end we divide the discussion into three cases, then apply
formula (4.7).

i) For $(G,p)=(G_{2},2),(F_{4},2),(E_{6},2),(F_{4},3),(E_{6},3),(E_{7},3)$ or
$(E_{8},5)$ the set $D_{1}(G,p)$ is a singleton $\{t_{i}\}$. Thus, according
to (4.7), for each $(G,p)$ there is just one relation $\mathcal{H}_{i,K}$
which takes the form $\rho_{l(g_{i})-1}\cdot x_{t_{i}}=0$. Precisely, these
relations are

\begin{center}%
\begin{tabular}
[c]{l||l|l|l|l|l|l|l}\hline
$(G,p)$ & $(G_{2},2)$ & $(F_{4},2)$ & $(E_{6},2)$ & $(F_{4},3)$ & $(E_{6},3)$
& $(E_{7},3)$ & $(E_{8},5)$\\\hline
$\mathcal{H}_{i,K}$ & $x_{6}\varrho_{11}$ & $x_{6}\varrho_{11}$ &
$x_{6}\varrho_{11}$ & $x_{8}\varrho_{23}$ & $x_{8}\varrho_{23}$ &
$x_{8}\varrho_{23}$ & $x_{12}\varrho_{59}$\\\hline
\end{tabular}
.
\end{center}

\noindent These have been stated in Theorem B.

ii) For $(G,p)=(E_{8},3)$ we have $D_{1}(G,p)=\{8,20\}$. By (4.7) there are
$6$ relations of the type $\mathcal{H}_{i,K}$, which are presented below:

\begin{center}
{\footnotesize \renewcommand{\arraystretch}{0.8}
\begin{tabular}
[c]{l||l|l|l|l|l|l}\hline
$(i,K)$ & $(8,\{8\})$ & $(8,\{20\})$ & $(8,\{8,20\})$ & $(20,\{8\})$ &
$(20,\{20\})$ & $(20,\{8,20\})$\\\hline
$\mathcal{H}_{i,K}$ & $\rho_{23}\cdot x_{8}$ & $\rho_{23}\cdot x_{20}%
-x_{8}^{2}\cdot\mathcal{C}_{\left\{  8,20\right\}  }$ & $\rho_{23}%
\cdot\mathcal{C}_{\left\{  8,20\right\}  }$ & $\rho_{59}\cdot x_{8}-x_{20}%
^{2}\cdot\mathcal{C}_{\left\{  8,20\right\}  }$ & $\rho_{59}\cdot x_{20}$ &
$\rho_{59}\cdot\mathcal{C}_{\left\{  8,20\right\}  }$\\\hline
\end{tabular}
,}
\end{center}

\noindent These have been recorded in v) of Theorem B.

iii) For the remaining cases $(G,p)=(E_{7},2)$ or $(E_{8},2)$ we have
$D_{1}(E_{7},2)=\{6,10,18\}$ and $D_{1}(E_{8},2)=\{6,10,18,30\}$,
respectively. Since the relations $\mathcal{H}_{i,K}$ are too many, we have
simply referred them, in iv) and v) of Theorem B, to the general formula (4.7).

The proof of theorem B has now been completed.\hfill$\square$

\bigskip

\noindent\textbf{Remark 4.9.} Granted with the general results in Theorems
4.1, 4.4, 4.7 and 4.9, one can recover the computations \cite{B1,P} of the
cohomology $H^{\ast}(G)$ for the classical groups $G=SU(n)$, $Spin(n)$ or
$Sp(n)$. For instance, if $G=SU(n)$ or $Sp(n)$ one gets from $D_{1}%
(G,p)=\emptyset$ (by Table 2) and Theorem 4.7 that $\tau_{p}(G)=0$ for any
prime $p$. It implies by Theorem 4.4 that

\begin{quote}
$H^{\ast}(G)=\Delta(\rho_{2l_{1}-1},\cdots,\rho_{2l_{n}-1})$,
\end{quote}

\noindent where $\rho_{2l-1}^{2}\in\tau_{2}(G)=0$ by Lemma 2.2. Thus,
inputting the value of the degree set $q(G)=\left\{  l_{1},\cdots
,l_{n}\right\}  $ for $G=SU(n)$ or $Sp(n)$, one gets the following results of
Borel \cite{B1}

\begin{quote}
$H^{\ast}(SU(n))=\Lambda(\rho_{3},\cdots,\rho_{2n-1})$, $H^{\ast
}(Sp(n))=\Lambda(\rho_{3},\cdots,\rho_{4n-1})$.\hfill$\square$
\end{quote}

\section{Historical remarks}

\textbf{Remark 5.1. }The study of the cohomology of Lie groups began with
coefficients over a field $\mathbb{F}$. Notably, Hopf \cite{H} observed, for
the first time, that the group product $\mu:G\times G\rightarrow G$ on $G$
induces a map of algebras

\begin{enumerate}
\item[(1.2)] $\mu^{\ast}:H^{\ast}(G;\mathbb{F})\rightarrow H^{\ast}(G\times
G;\mathbb{F})\cong H^{\ast}(G;\mathbb{F})\otimes H^{\ast}(G;\mathbb{F})$,
\end{enumerate}

\noindent where the isomorphism follows from the K\"{u}nneth formula. It
furnishes the cohomology $H^{\ast}(G;\mathbb{F})$ with an additional
co-produce $\mu^{\ast}$, making the pair $\left\{  H^{\ast}(G;\mathbb{F}%
),\mu^{\ast}\right\}  $ an \textsl{Hopf algebra} over $\mathbb{F}$. In term of
$\mu^{\ast}$ one can specify the set of so called "\textsl{primitive
elements}"\textsl{ }of\textsl{ }$H^{\ast}(G;\mathbb{F})$

\begin{quote}
$P(G;\mathbb{F})=\{a\in H^{\ast}(G;\mathbb{F})\mid\mu^{\ast}(a)=a\otimes1$
$\oplus1\otimes a\}$,
\end{quote}

\noindent which is clearly a linear subspace of $H^{\ast}(G;\mathbb{F})$. For
the case $\mathbb{F}$ is the field $\mathbb{R}$ of reals Hopf \cite{H} proved that

\bigskip

\noindent\textbf{Theorem 5.1. }\textsl{If }$\left\{  z_{1},\cdots
,z_{k}\right\}  $\textsl{ is a homogeneous basis of }$P(G;\mathbb{R}%
)$\textsl{, then}

\begin{enumerate}
\item[(5.2)] $H^{\ast}(G;\mathbb{R})=\Lambda_{\mathbb{R}}(z_{1},\cdots,z_{k}%
)$\textsl{, }$\deg z_{i}\equiv1\operatorname{mod}2$.\hfill$\square$
\end{enumerate}

Borel \cite{B1} initiated the investigation of the algebra $H^{\ast
}(G;\mathbb{F}_{p})$ over a finite field $\mathbb{F}_{p}$. Following the idea
of Hopf he started with a classification on the Hopf algebra $H^{\ast
}(G;\mathbb{F}_{p})$ over a finite field $\mathbb{F}_{p}$.

\bigskip

\noindent\textbf{Theorem} \textbf{5.2}. \textsl{If }$\left\{  z_{1}%
,\cdots,z_{k}\right\}  $\textsl{ is a homogeneous basis of} $P(G;\mathbb{F}%
_{p})$\textsl{, then}

\begin{enumerate}
\item[(5.3)] $H^{\ast}(G;\mathbb{F}_{p})=B(z_{1})\otimes\cdots\otimes
B(z_{k})$\textsl{,}
\end{enumerate}

\noindent\textsl{where each factor }$B(z_{i})$\textsl{ is one of the next
"monogenic Hopf algebra" over }$\mathbb{F}_{p}$\textsl{:}

\begin{quote}%
\begin{tabular}
[c]{l|l|l}\hline\hline
$B(z)$ & $p=2$ & $p\neq2$\\\hline
$\deg(z)=even$ & $\mathbb{F}_{2}(x)/\left\langle x^{2^{r}}\right\rangle $ &
$\mathbb{F}_{p}(x)/\left\langle x^{p^{r}}\right\rangle $\\\hline
$\deg(z)=odd$ & $\mathbb{F}_{2}(x)/\left\langle x^{2^{r}}\right\rangle $ &
$\Lambda_{\mathbb{F}_{p}}(x)$\\\hline\hline
\end{tabular}
.\hfill$\square$
\end{quote}

Based on Theorems 5.1 and 5.2 the algebras $H^{\ast}(G;\mathbb{F})$ were
largely computed by Borel \cite{B2,B3}, Chevalley \cite{Ch0}, Araki
\cite{AS,A1,A2,A3}, Toda, Kono, Mimura and Shimada \cite{T,KM1,KM2,KM3,KMS,Ko}%
, case by case.

Schubert calculus is an inspiring subject. It promoted the enumerative
geometry of the 19th century growing into the algebraic geometry founded by
Van der Waerden and Andr\'{e} Weil, and has integrated into many branches of
mathematics, such as the theory of characteristic classes \cite{MS}, the
string theory \cite{Ka}, and algebraic combinatorics \cite{Fu}. The present
work, together with the earlier one \cite{DZ1}, extends the calculus to a
unified construction of the cohomologies of all compact and simply-connected
Lie groups.

\bigskip

\textbf{Remark 5.2.} A compact Lie group $K$ is called \textsl{semi-simple} if
it is isomorphic to the quotient group $G/L$, where $G$ is compact and
simply-connected, and $L\subset G$ is a subgroup of the center of $G$. In
particular, the orthogonal group $SO(n)$, the projective unitary group
$PU(n):=SU(n)/\mathbb{Z}_{n}$, the adjoint exceptional Lie groups
$Ad(E_{6})=E_{6}/\mathbb{Z}_{3}$ and $Ad(E_{7})=E_{7}/\mathbb{Z}_{2}$, are
example non simply-connected semi-simple Lie groups, e.g. \cite{BB}.

To simplify the presentation we have assumed, at the beginning, that the Lie
groups $G$ under our consideration are the compact and simply-connected ones.
We mention at this point that the method developed in this paper applies
equally well to construct the cohomology of all compact and semi-simple Lie
groups $K$. Precisely, let $T^{\prime}$ be a maximal torus of $K$, and let $T$
be a maximal torus of $G$ that corresponds to $T^{\prime}$ under the universal
covering $G\rightarrow K$. Then, in view of the canonical isomorphism
$K/T^{\prime}=G/T$ in flag manifolds, the second page of the Serre spectral
sequence $\{E_{r}^{\ast,\ast}(K),d_{r}^{\prime}\}$ of $K\rightarrow
K/T^{\prime}$ is

\begin{quote}
$E_{2}^{\ast,\ast}(K)=H^{\ast}(G/T)\otimes H^{\ast}(T^{\prime})$,
\end{quote}

\noindent where

i) the Schubert presentation of the ring $H^{\ast}(G/T)$ is available by (1.7);

ii) the differential $d_{2}^{\prime}:H^{1}(T^{\prime})\rightarrow H^{2}(G/T)$
can be formulated in term of the Cartan matrix of $K$, and the fundamental
group $\pi_{1}(K)=L$ (e.g. \cite[Theorem 2.3]{D1}).

\noindent Combining these ideas the cohomologies $H^{\ast}(K)$ of all
semi-simple Lie groups can been constructed using the methods developed in
this paper (e.g. \cite{D2,D3}).

For extension of Schubert calculus to computing with the integral cohomology
of certain homogeneous spaces, we refer to \cite[Section \S 5]{DZ0}.

\bigskip

\textbf{Remark 5.3.} The approach to the cohomology of Lie groups, using the
Serre spectral sequence of $\pi:G\rightarrow G/T$, has been interested by
several authors.

In \cite{L} Leray proved that the Serre spectral sequence of $\pi$ satisfies

\begin{quote}
$E_{3}^{\ast,\ast}(G;\mathbb{R})=E_{\infty}^{\ast,\ast}(G;\mathbb{R})$.
\end{quote}

\noindent Later on, Ka\v{c} and Marlin \cite{K,M2} conjectured that there
exist additive isomorphisms

\begin{enumerate}
\item[(5.4)] $E_{3}^{\ast,\ast}(G)\cong E_{\infty}^{\ast,\ast}(G)$ and
$E_{\infty}^{\ast,\ast}(G)\cong H^{\ast}(G)$.
\end{enumerate}

\noindent Since our construction of $H^{\ast}(G)$ factor through $E_{3}%
^{\ast,\ast}(G)$ (see in Section \S 1.2), these conjectures can be confirmed
within our context.

Combining Schubert calculus on $G/T$, the invariant theory of Lie groups
\cite{Ch1}, and the Serre spectral sequence of $\pi$ with real coefficients,
Reeder \cite{R} gave a complete treatment of the following computation of
Borel, Chevalley and Leray \cite{B1,Ch0,L}

\begin{enumerate}
\item[(5.5)] $H^{\ast}(G;\mathbb{R})=\Lambda_{\mathbb{R}}(\xi_{2l_{1}%
-1},\cdots,\xi_{2l_{n}-1})$, $q(G)=\{l_{1},\cdots,l_{n}\}$.
\end{enumerate}

\noindent In our context, taking $\xi_{2l-1}=\rho_{2l-1}\otimes1\in H^{\ast
}(G)\otimes\mathbb{R}$, formula (5.5) follows from Theorem 4.1.

\bigskip

\begin{quote}
Haibao Duan,

dhb@math.ac.cn,

Yau Mathematical Science Center, Tsinghua University, Beijing 100084;

Academy of Mathematics and Systems Sciences, Chinese Academy of Sciences,
Beijing 100190.
\end{quote}

\end{document}